\documentclass[11pt]{article}
\usepackage{graphicx}
\usepackage{color}
\usepackage{amsmath,amssymb,amsthm,amsfonts}
\usepackage{amssymb}
\newtheorem{thm}{Theorem}[section]

\newtheorem{lem}[thm]{Lemma}
\newtheorem{prop}[thm]{Proposition}

\theoremstyle{definition}

\theoremstyle{remark}

\numberwithin{equation}{section}

\DeclareMathSymbol{\C}{\mathalpha}{AMSb}{"43}

\newcommand{\R}{{\mathbb{R}}}

\newcommand{\bsub}{\begin{subequations}}
\newcommand{\esub}{\end{subequations}$\!$}

\begin{document}
\title{Complex-valued Solutions of the Planar Schr\"odinger-Newton System }
\author{Hongfei Zhang\thanks{School of Mathematics and Statistics, Central China Normal University, P.O. Box 71010, Wuhan 430079,
P. R. China.  Email: \texttt{hfzhang@mails.ccnu.edu.cn}.}
 \, and\, Shu Zhang\thanks{School of Mathematics and Statistics, Central China Normal University, P.O. Box 71010, Wuhan 430079,
P. R. China.  Email: \texttt{shu@mails.ccnu.edu.cn}.}
}
\date{\today}

\smallbreak \maketitle
\begin{abstract}
In this paper, we consider complex-valued solutions of the planar Schr\"odinger-Newton system, which can be described by minimizers of the constraint minimization
problem. It is shown that there exists a critical rotational velocity $0<\Omega^*\leq \infty$, depending on the general trapping potential $V(x)$, such that for any rotational velocity $0\leq\Omega<\Omega^*$, minimizers exist if and only if $0<a<a^*:=\|Q\|_{2}^{2}$, where $Q>0$
is the unique positive solution of $-\Delta u+u-u^3=0$ in $\R^2$. Moreover, under some suitable assumptions on $V(x)$, applying blow-up analysis and energy estimates, we present a detailed analysis on the concentration behavior of minimizers as $a\nearrow a^*$.
\end{abstract}

\section{Introduction}
In this paper, we study complex-valued solutions of the following Schr\"odinger-Newton system
\begin{equation}\label{0.1}
\left\{
\begin{aligned}
i\psi_{t} &=-\Delta\psi+V(x)\psi+2\pi \phi\psi-|\psi|^{2}\psi\ \ \hbox{in} \ \ \R^2\times\R,\\
\Delta \phi&=|\psi|^2\ \ \hbox{in}\ \ \R^2\times\R,
\end{aligned}
\right.
\end{equation}
where $\psi:\R^2\times\R\rightarrow\C$ is the time-dependent wave function,
$x\mapsto V(x)$ is a real external potential. The function $\phi$ represents an internal potential for a nonlocal self-interaction of the wave function $\psi$, and the nonlinear term $|\psi|^{2}\psi$ is frequently used to model the interaction among particles (cf. \cite{Silvia1}).
The system \eqref{0.1} arises in several areas, such as quantum mechanics and semiconductor physics, we refer the reader to  \cite{Ando,Doniach,Fisher,Frohlich,Garnier,Markowich,Penrose1,Penrose2} for more details on its physical aspects.

Over the past few decades, the system \eqref{0.1} has attracted considerable attention due to its physical relevance (cf. \cite{CaoD,Cassani,Choquard,Cingolani,Silvia1,GLL,HMT,Masaki,Stubbe,SV,Wang}). More specifically, taking no account of the nonlinear term $|\psi|^{2}\psi$ and the external potential $V(x)$ in (\ref{0.1}), Harrison, Moroz and Tod \cite{HMT} established the existence of bound states for (\ref{0.1}) by numerical study, after which, applying a shooting method,
Choquard, Stubbe and Vuffray \cite{Choquard} proved the existence of a unique positive radially symmetric solution for (\ref{0.1}). Recently, for the standing wave ansatz $\psi(x,t)=e^{i\mu t}u(x)$, $\mu\in\mathbb{R}$,
Cingolani and Weth \cite{Silvia1} deduced a strong compactness condition for Cerami sequences, based on which they established various existence results of
solutions for (\ref{0.1}) with $V(x)\in L^{\infty}(\mathbb{R}^{2})$.
Moreover, employing
the variational method, Guo, Liang and Li \cite{GLL} investigated the existence and the correlation properties of solutions for
(\ref{0.1}) in the case of $V(x)=\ln(1+|x|^2)$.

However, there are few researches for complex-valued solutions of (\ref{0.1}), and the study of complex-valued solutions can help us better understand the properties of solutions (cf. \cite{HMT,SV}). Inspired by the above facts, we shall seek solutions with no
counterpart among those found already by making an ansatz of rigid rotation.
That is, we attempt to study solutions of (\ref{0.1}) which in polar coordinates $r, \theta$ take the form
\begin{equation}\label{1.3}
\psi(r,\theta,t)=e^{-i\mu t}u(r,\theta+\Omega t)
\end{equation}
for real constants $\mu$ and $\Omega$. To separate the Schr\"odinger-Newton system (\ref{0.1}) for a solution like (\ref{1.3}), we go into the rotating frame with Cartesian coordinates $y_{1}, y_{2}$ given by
\begin{equation*}
\begin{pmatrix}
y_1\\
y_2
\end{pmatrix}
=
\begin{pmatrix}
cos \Omega t&sin\Omega t\\
-sin\Omega t&cos\Omega t
\end{pmatrix}
\begin{pmatrix}
x_1\\
x_2
\end{pmatrix}
:=A
\begin{pmatrix}
x_1\\
x_2
\end{pmatrix},
\end{equation*}
and rename the potential $V(A^{-1}y)$ to be $V(y)$, then
the time-dependence separates off to leave the system \eqref{0.1} as
\begin{equation}\label{1.1}
\begin{aligned}
&-\Delta u+V(y)u+\Big(\int_{\mathbb{R}^{2}}\ln|y-z||u(z)|^2dz\Big)u+
i\Omega(y^{\bot}\cdot\nabla u)\\
&=\mu u+|u|^{2}u\ \ \text{in}\ \ \mathbb{R}^{2},\ \ \text{where}\ \ y^{\perp}=(-y_{2},y_{1}).
\end{aligned}
\end{equation}

It is well known that solutions of (\ref{1.1}) can be obtained as minimizers of the following  complex-valued variational problem
\begin{equation}\label{1.4}
e(a):=\inf_{u\in S_{a}}E_{a}(u),\ \ a>0,
\end{equation}
where the energy functional $E_{a}(u)$ is defined by
\begin{equation}\label{1.5}
\begin{split}
E_{a}(u):=&\int_{\R^2}(|\nabla u|^2+V(x)|u|^2)dx+\frac{1}{2}\int_{\R^2}\int_{\R^2}\ln|x-y||u(x)|^2|u(y)|^2dxdy\\
&-\frac{1}{2}\int_{\R^2}|u(x)|^4dx
-\Omega\int_{\R^2}x^{\bot}\cdot(iu, \nabla u)dx,
\end{split}
\end{equation}
here $(iu,\nabla u)=\frac{i(u\nabla\bar{u}-\bar{u}\nabla u)}{2}$.
Stimulated by \cite{Silvia1,GLL,Stubbe}, we consider the space $S_a$ of $e(a)$ satisfying
\begin{equation*}
S_a:=\big\{u\in \mathcal{H}: \|u\|^2_2=a\big\},\ \ \mathcal{H}=\Big\{u\in H^1(\R^2, \C):\int_{\R^2}V(x)|u(x)|^2dx<\infty\Big\},
\end{equation*}
with the associated norm
\begin{equation}\label{1.9}
\|u\|_\mathcal{H}=\left\{\int_{\R^2}\left[|\nabla u|^2+\big(1+V(x)\big)|u|^2\right]dx\right\}^{\frac{1}{2}}.
\end{equation}

In recent years, there are many results (cf. \cite{Dinh,Guo,GLL,GLP,GLY}) related to the variational problem (\ref{1.4}). In particular, when there is no logarithmic convolution term  in (\ref{1.4}), Guo and his collaborators \cite{Guo,GLP,GLY} studied various quantitative properties of minimizers for (\ref{1.4}). In addition, Guo, Liang and Li \cite{GLL} studied the existence and uniqueness of minimizers for (\ref{1.4}) with $\Omega=0$ and $V(x)=\ln(1+|x|^2)$. Compared with these existing results, it is more interesting and challenging when the logarithmic convolution term and the rotation term are simultaneously involved. We need to carefully analyze the dominant role and mutual influence of these two items. On the one hand, it should be noted that the Newton potential in two-dimensional space is sign-changing, and it is singular in points at origin and infinity. On the other hand, minimizers are complex-valued when considering the rotation term, which makes the calculation more difficult.


In this paper, we shall study the variational problem (\ref{1.4}) to obtain the existence and related properties of complex-valued solutions for the Schr\"odinger-Newton system (\ref{0.1}).  Before stating our main results, we first recall the following Gagliardo-Nirenberg inequality (cf. \cite{Weinstein})
\begin{equation}\label{1.13}
\int_{\R^2}|u(x)|^4dx\leqslant\frac 2{\|Q\|_2^2}\int_{\R^2}|\nabla u(x)|^2dx\int_{\R^2}|u(x)|^2dx,\ \ u\in H^1(\R^2,\R),
\end{equation}
and the equality in (\ref{1.13}) is achieved at $u(x)=mQ(mx)$ ($m\neq 0$ is arbitrary), where $Q=Q(|x|)$ is the unique (up to translations) positive radial solution of the following nonlinear scalar field equation (cf. \cite{Gidas,Kwong})
\begin{equation}\label{1.6}
-\Delta u+u-u^3=0\ \ {\rm in}\ \ \R^2,\ \ {\rm where}\ \ u\in H^1(\R^2,\R).
\end{equation}
Note also from \cite[Proposition 4.1]{Gidas} that $Q=Q(|x|)>0$ strictly decreases in $|x|$ and admits the following exponential decay
\begin{equation}\label{1.15}
Q(x),~|\nabla Q(x)|=O(|x|^{-\frac{1}{2}}e^{-|x|})\ \ \text{as}\ \ |x|\rightarrow\infty.
\end{equation}
Moreover, it follows from \cite[Lemma 8.1.2]{Cazenave} that $Q$ satisfies
\begin{equation}\label{1.16}
\int_{\mathbb{R}^{2}}|\nabla Q|^{2}dx=\int_{\mathbb{R}^{2}}Q^{2}dx
=\frac{1}{2}\int_{\mathbb{R}^{2}}Q^{4}dx.
\end{equation}

Throughout this paper, we always assume that the external potential $V(x)$ satisfies
\begin{equation}\label{1.10}
0\leq V(x)\in L_{loc}^{\infty}(\mathbb{R}^{2}),\ \ \lim\limits_{|x|\rightarrow\infty}\frac{V(x)}{|x|^{2}}>0,
\end{equation}
and we define the critical rotational velocity $\Omega^{*}$ by
\begin{equation}\label{1.11}
\Omega^{*}:=\sup\Big\{\Omega>0:V(x)-\frac{\Omega^{2}}{4}|x|^{2}\rightarrow\infty \ \ \text{as} \ \ |x|\rightarrow\infty\Big\}.
\end{equation}
Using the above notations,
we have the following result concerning the existence and nonexistence of minimizers for $e(a)$. 
\begin{thm}\label{thm:1.1}
Let $Q$ be the unique positive radial solution of (\ref{1.6}). Suppose $V(x)$ satisfies \eqref{1.10}, then we have
\begin{enumerate}
\item If $0\leq\Omega<\Omega^*$ and $0<a<a^*:=\|Q\|_{2}^{2}$, there exists at least one minimizer for $e(a)$;

\item If $0\leq\Omega<\Omega^*$ and $a\geq a^*:=\|Q\|_{2}^{2}$, there is no minimizer for $e(a)$;

\item  If $\Omega>\Omega^*$, then for any $a>0$, there is no minimizer for $e(a)$.
\end{enumerate}
Moreover, we have  $\lim\limits_{a\nearrow a^*}e(a)=-\infty$ for any $\Omega\geq0$.
\end{thm}
From the above conclusion, we see that the existence and
nonexistence of minimizers for $e(a)$ are the same as the case  without rotation, provided that the rotating velocity $\Omega$ is smaller than the critical value $\Omega^*$.
Note that the appearance of the rotation term increases the difficulty of proof. We mainly apply the following diamagnetic inequality (cf.\cite{Lieb1}) to overcome this obstacle. For $\mathcal{A}=\frac{\Omega}{2}x^{\bot}$,
\begin{equation}\label{1.12}
|\nabla u|^{2}-\Omega x^{\bot}\cdot(iu, \nabla u)+\frac{\Omega^{2}}{4}|x|^{2}|u|^{2}=|(\nabla-i\mathcal{A})u|^{2}\geq\Big|\nabla|u|\Big|^{2},\ \ \ u\in H^{1}(\mathbb{R}^{2}, \mathbb{C}).
\end{equation}
Moreover, our analysis shows that the proof of Theorem \ref{thm:1.1} needs the
following Hardy-Littlewood-Sobolev inequality (cf. \cite{Lieb1})
\begin{equation}\label{1.17}
\int_{\R^2}\int_{\R^2}\frac{1}{|x-y|}|u(x)||v(y)|dxdy\leqslant C\|u\|_{\frac{4}{3}}\|v\|_{\frac{4}{3}},\ \ \hbox{where}\ \ u,\ v\in L^{\frac{4}{3}}(\R^2).
\end{equation}
Actually, the inequality \eqref{1.17} is often used throughout this paper to handle with the logarithmic convolution term of $E_{a}(u)$.

By the variational argument, if $u_a$ is a minimizer of $e(a)$ for $a\in(0, a^*)$, then there exists a Lagrange multiplier $\mu_a$ satisfies the following Euler-Lagrange equation
\begin{equation}\label{1.18}
\begin{split}
&-\Delta u_a+V(x)u_a+\Big(\int_{\mathbb{R}^{2}}\ln|x-y||u_a(y)|^2dy\Big)u_a+
i\Omega(x^{\bot}\cdot\nabla u_a)\\
&=\mu_{a}u_a+|u_a|^{2}u_a\ \ \text{in}\ \ \mathbb{R}^{2},
\end{split}
\end{equation}
i.e. $u_a$ is a solution of (\ref{1.1}) with $\mu=\mu_a$. By making full use of the equation (\ref{1.18}), we next focus on investigating the concentration behavior of minimizers for $e(a)$ as $a\nearrow a^*$. In addition to the assumption that $V(x)$ satisfies (\ref{1.10}), we further assume that $V_{\Omega}(x)=V(x)-\frac{\Omega^{2}}{4}|x|^2$ satisfies the following conditions:
\begin{enumerate}
\item [\rm $(V)$.]
$0\leq V_{\Omega}(x)\in C_{loc}^{\alpha}(\mathbb{R}^{2})~(0<\alpha<1)$, $\{x\in \mathbb{R}^2: V_{\Omega}(x)=0\}=\{0\}$ and $V(x)\leq Ce^{\gamma|x|}$ for some $\gamma>0$ as $|x|\rightarrow\infty$.
\end{enumerate}
Under the assumptions (\ref{1.10}) and $(V)$, the main result of concentration behavior of minimizers for $e(a)$ is stated as follows.
\begin{thm}\label{thm:1.2}
Suppose that $V(x)$ satisfies (\ref{1.10}) and $(V)$, and assume $0<\Omega<\Omega^*$, where $\Omega^*>0$ is defined as in (\ref{1.11}). If $u_a$ is a minimizer of $e(a)$, then we have
\begin{equation}\label{1.21}
\begin{aligned}
w_a(x):=\frac{2}{a^*}\Big(\frac{a^*-a}{a^*}\Big)^{\frac{1}{2}}u_a\Big(\frac{2}{a^*}
\big(\frac{a^*-a}{a^*}\big)^{\frac{1}{2}}x+x_a\Big)
e^{-i\big(\frac{\Omega}{a^*}\big(\frac{a^*-a}{a^*}\big)^{\frac{1}{2}}x\cdot x_{a}^{\bot}-\theta_{a}\big)}\\
\rightarrow\frac{1}{\sqrt{a^*}}Q\Big(\frac{1}{\sqrt{a^*}}x\Big)\ \ \text{as}\ \ a\nearrow a^*
\end{aligned}
\end{equation}
strongly in $H^1(\mathbb{R}^2,\mathbb{C})$ $\cap$ $L^\infty(\mathbb{R}^2,\mathbb{C})$, where $\theta_{a}\in[0,2\pi)$ is a properly chosen constant. Here $x_{a}\in\mathbb{R}^{2}$ is the unique global maximal point of $|u_a|$ and $\lim_{a\nearrow a^*}V_{\Omega}(x_a)=0$.
\end{thm}

Theorem \ref{thm:1.2} gives a detailed description of the concentration behavior of minimizers for $e(a)$. The proof of Theorem \ref{thm:1.2} relies heavily on the refined energy estimates of $e(a)$. Since both the rotation term and the logarithmic convolution term appear in $e(a)$, there exist some extra difficulties in the proof of Theorem \ref{thm:1.2}. Firstly, note that $u_a$ is a complex-valued function, so we cannot use the Gagliardo-Nirenberg inequality (\ref{1.13}) directly to $u_a$ as in \cite{GLL}. For this reason, stimulated by \cite{GLL0}, we shall first establish the asymptotic behavior of $|u_a|$ as $a\nearrow a^*$. Secondly, owing to $\lim_{a\nearrow a^*}e(a)=-\infty$, compared with \cite{GLL0, GLY}, in Section 3 we use different analytical methods to derive the uniform boundedness of the scaled minimizers $w_{a}$ as $a\nearrow a^*$. Finally, since $e(a)$ has both the rotation term and the logarithmic convolution term, we need to explore new techniques to get the uniqueness and estimates of the global maximum point for $|u_a|$.

This paper is organized as follows. In Section 2, we first give some preliminary results, after which we prove Theorem \ref{thm:1.1} on the existence and non-existence of minimizers for $e(a)$. In Section 3, we shall prove Theorem \ref{thm:1.2} on the concentration behavior of minimizers for $e(a)$ by employing blow-up analysis and energy estimates. We finally give in Appendix the detailed proof of some results used in the proof of Lemma \ref{pro:3.2}.

\section{Existence of minimizers for $e(a)$}
In this section, we first establish various preliminary results which are often used throughout the whole paper. We then follow them to finish the proof of Theorem \ref{thm:1.1} in subsection 2.2.

\subsection{Preliminaries}

Firstly, we introduce the following compactness lemma.
\begin{lem}\label{lem:2.1}
Suppose that $V(x)\in L_{loc}^{\infty}(\mathbb{R}^{2})$ and $\lim\limits_{|x|\rightarrow\infty}V (x)=\infty$. If $2\leq q<\infty$, then the embedding $\mathcal{H}\hookrightarrow L^{q}(\mathbb{R}^{2}, \mathbb{C})$ is compact.
\end{lem}
\noindent The proof of Lemma \ref{lem:2.1} is similar to those in \cite{ZhangJ} and the references therein, so we
omit it here.

Next, as performed in \cite{Silvia1}, we denote the space by
$$X:=\{u\in H^{1}(\mathbb{R}^{2}): \int_{\mathbb{R}^{2}}\ln\big(1+|x|\big)|u(x)|^{2}dx<\infty\}.$$
By the assumption of (\ref{1.10}), we can obtain that $|u|\in X$ if $u\in\mathcal{H}$.
We then define the following symmetric bilinear forms
\begin{equation*}
B_1(u,v):=\int_{\R^2}\int_{\R^2}{\rm ln}\big(1+|x-y|\big)u(x)v(y)dxdy,
\end{equation*}
\begin{equation*}
B_2(u,v):=\int_{\R^2}\int_{\R^2}{\rm ln}\Big(1+\frac1{|x-y|}\Big)u(x)v(y)dxdy,
\end{equation*}
and
\begin{equation}
\begin{aligned}\label{200}
B_0(u,v):=B_1(u,v)-B_2(u,v)=\int_{\R^2}\int_{\R^2}{\rm ln}|x-y|u(x)v(y)dxdy,\ \ u,v\in X.
\end{aligned}
\end{equation}
Since
\begin{equation*}
{\rm ln}\big(1+|x-y|\big)\leqslant{\rm ln}\big(1+|x|+|y|\big)\leqslant{\rm ln}\big(1+|x|\big)+{\rm ln}\big(1+|y|\big),\ \ x,y\in\R^2,
\end{equation*}
we have the estimate
\begin{equation}
\begin{aligned}\label{2.1}
\left|B_1(uv,wz)\right|\leqslant&\int_{\R^2}\int_{\R^2}\Big[{\rm ln}(1+|x|)+{\rm ln}(1+|y|)\Big]\big|u(x)v(x)w(y)z(y)\big|dxdy\\
\leqslant&\|u\|_*\|v\|_*\|w\|_2\|z\|_2+\|u\|_2\|v\|_2\|w\|_*\|z\|_*,\ \ u,v,w,z\in X,
\end{aligned}
\end{equation}
where
\begin{equation*}
\|u\|_*:=\left\{\int_{\R^2}{\rm ln}(1+|x|)|u(x)|^2dx\right\}^{\frac{1}{2}},\ \ u\in X.
\end{equation*}
Note that $0<{\rm ln}(1+r)< r$ for any $r>0$, it follows from the Hardy-Littlewood-Sobolev inequality \eqref{1.17} that there exists a constant $C>0$ such that
\begin{equation}
\begin{aligned}\label{2.2}
|B_2(u,v)|\leqslant&\int_{\R^2}\int_{\R^2}{\rm ln}\Big(1+\frac1{|x-y|}\Big)|u(x)||v(y)|dxdy\\
\leqslant&\int_{\R^2}\int_{\R^2}\frac1{|x-y|}|u(x)||v(y)|dxdy\\
\leqslant&C\|u\|_{\frac{4}{3}}\|v\|_{\frac{4}{3}},\ \ u,v\in L^{\frac{4}{3}}(\R^2).
\end{aligned}
\end{equation}

Applying the compactness Lemma \ref{lem:2.1},  we derive from  \eqref{2.1} and \eqref{2.2} that the functionals $B_i(u^2,v^2)$ are well defined on $X\times X$ for $i=0,1,2$. Moreover, denote functionals on $X$ as follows
\begin{equation*}
V_1(u):=B_1(u^2,u^2)=\int_{\R^2}\int_{\R^2}{\rm ln}\big(1+|x-y|\big)u^2(x)u^2(y)dxdy,
\end{equation*}
 \begin{equation*}
V_2(u):=B_2(u^2,u^2)=\int_{\R^2}\int_{\R^2}{\rm ln}\Big(1+\frac1{|x-y|}\Big)u^2(x)u^2(y)dxdy,
\end{equation*}
and
\begin{equation*}
V_0(u):=B_0(u^2,u^2)=\int_{\R^2}\int_{\R^2}{\rm ln}|x-y|u^2(x)u^2(y)dxdy.
\end{equation*}
Combining (\ref{200})--(\ref{2.2}), we obtain that
\begin{equation}
\begin{aligned}\label{22}
&\int_{\R^2}\int_{\R^2}{\rm ln}|x-y|u^2(x)u^2(y)dxdy\\
\leq&\int_{\R^2}\int_{\R^2}{\rm ln}\big(1+|x-y|\big)u^2(x)u^2(y)dxdy\\
&+
\int_{\R^2}\int_{\R^2}{\rm ln}\Big(1+\frac1{|x-y|}\Big)u^2(x)u^2(y)dxdy\\
\leq&2\|u\|_*^{2}\|u\|_2^{2}+C\|u\|_{\frac{8}{3}}^{4},\ \ u\in X.
\end{aligned}
\end{equation}
Using the above notations, recall from \cite[Lemma 2.2]{Silvia1} that $V_i(\cdot)~(i=0,1,2)$ satisfy the following crucial properties.

\begin{lem}\label{lem:2.2} The space $X$ and the functionals $V_i(\cdot)$ $(i=0,1,2)$ have the following properties

\begin{enumerate}
		\item  The space $X$ is compactly embedded into $L^s(\R^2)$ for all $s\in[2,\infty)$;

\item The functionals $V_0(\cdot), V_1 (\cdot)$ and $V_2(\cdot)$ are of class $C^1$ on $X$, and
 $$\langle V_i'(u), v\rangle=4B_i(u^2, uv)\ \  \mbox{for}\ \ u, v\in X, \ \
  i=0, 1, 2 ;$$

\item $V_2(\cdot)$ is continuously differentiable on $L^{\frac{8}{3}}(\R^2)$;

\item $V_1(\cdot)$ is weakly lower semicontinuous on $H^1(\R^2)$;

\item $V_0(\cdot)$ is weakly lower semicontinuous on $X$.
\end{enumerate}
\end{lem}

\subsection{The proof of Theorem \ref{thm:1.1}}

In this subsection, we shall prove Theorem \ref{thm:1.1} on the existence and nonexistence of minimizers for $e(a)$.
First of all, making full us of the diamagnetic inequality (\ref{1.12}) and the compactness Lemma \ref{lem:2.1}, we can get the following existence.

\begin{prop}\label{pro2.4}
Let $Q=Q(|x|)$ be the unique positive solution of (\ref{1.6}). Suppose $V(x)$ satisfies (\ref{1.10}), and $\Omega^*>0$ is defined as in (\ref{1.11}). Then for any $0\leq\Omega<\Omega^*$ and $0<a<a^*$, there exists at least one minimizer for $e(a)$.
\end{prop}
\noindent\textbf{Proof.}
Suppose that $u\in S_a$, we then deduce from the Gagliardo-Nirenberg inequality (cf.\cite{Weinstein}) that for any $p>2$,
\begin{equation}\label{2.0}
\|u\|_p\leq C_{p}a^{\frac{1}{p}} \Big(\int_{\R^2}\big|\nabla |u|\big|^2dx\Big)^{\frac{p-2}{2p}}.
\end{equation}
By the Hardy-Littlewood-Sobolev inequality \eqref{1.17}, we derive from (\ref{2.0}) that there exists a constant $C>0$ such that
\begin{equation*}
\int_{\R^2}\int_{\R^2}\frac{|u(x)|^2|u(y)|^2}{|x-y|}dxdy\leqslant C\Big(\int_{\R^2}|u(x)|^{\frac{8}{3}}dx\Big)^{\frac{3}{2}}\leqslant Ca^{\frac{3}{2}}\Big(\int_{\R^2}|\nabla |u||^2dx\Big)^{\frac{1}{2}},
\end{equation*}
it then follows from above that
\begin{equation}\label{2.6}
\begin{aligned}
&\int_{\R^2}\int_{\R^2}{\rm ln}\Big(1+\frac1{|x-y|}\Big)|u(x)|^2|u(y)|^2dxdy\\
\leqslant&\int_{\R^2}\int_{\R^2}\frac1{|x-y|}|u(x)|^2|u(y)|^2dxdy\\
\leqslant&Ca^{\frac{3}{2}}\Big(\int_{\R^2}|\nabla |u||^2dx\Big)^{\frac{1}{2}}.
\end{aligned}
\end{equation}

On the other hand, for any $0\leq\Omega<\Omega^*$, we deduce from (\ref{1.10}) and (\ref{1.11}) that
\begin{equation}\label{2.3}
\lim_{|x|\rightarrow\infty}\Big[V(x)-\frac{\Omega^2}{4}|x|^2\Big]=\infty.
\end{equation}
Since $V(x)\in L_{loc}^{\infty}(\R^2)$, it follows from (\ref{2.3}) that there exists a certain constant $M>0$ such that $V(x)-\frac{\Omega^2}{4}|x|^2>-M$ holds in $\R^2$, which implies that
$\int_{\R^2}(V(x)-\frac{\Omega^2}{4}|x|^2)|u|^2dx>-aM$. Combining the above facts, we deduce from (\ref{1.13}) and (\ref{1.12}) that
\begin{equation}\label{2.4}
\begin{split}
E_{a}(u)=&\int_{\R^2}\big(|\nabla u|^2+V(x)|u|^2\big)dx+\frac{1}{2}\int_{\R^2}\int_{\R^2}\ln|x-y||u(x)|^2|u(y)|^2dxdy\\
&-\frac{1}{2}\int_{\R^2}|u(x)|^4dx
-\Omega\int_{\R^2}x^{\bot}\cdot(iu, \nabla u)dx\\
=&\int_{\R^2}|(\nabla-i\mathcal{A})u|^2dx+
\frac{1}{2}\int_{\R^2}\int_{\R^2}\ln|x-y||u(x)|^2|u(y)|^2dxdy\\
&-\frac{1}{2}\int_{\R^2}|u(x)|^4dx+\int_{\R^2}\Big(V(x)-\frac{\Omega^2}{4}|x|^2\Big)|u(x)|^2dx\\
\geq&\int_{\R^2}\big|\nabla|u|\big|^2dx-\frac{1}{2}\int_{\R^2}\int_{\R^2}\ln\Big(1+\frac{1}{|x-y|}\Big)
|u(x)|^2|u(y)|^2dxdy\\
&-\frac{1}{2}\int_{\R^2}|u(x)|^4dx-aM\\
\geq&\frac{a^*-a}{a^*}\int_{\R^2}\big|\nabla|u|\big|^2dx-\frac{Ca^{\frac{3}{2}}}{2}
\Big(\int_{\R^2}\big|\nabla|u|\big|^2dx\Big)^\frac{1}{2}-aM\\
\geq&\frac{a^*-a}{2a^*}\int_{\R^2}\big|\nabla|u|\big|^2dx-\frac{Ca^{3}a^*}{a^*-a}-aM,
\end{split}
\end{equation}
which implies that $E_{a}(u)$ is bounded from below.

Let $\{u_{n}\}\subset S_{a}$ be a minimizing sequence of $e(a)$, we see from (\ref{2.4}) that $\int_{\R^2}\big|\nabla|u_n|\big|^2dx$ is bounded uniformly in $n$. Applying the Gagliardo-Nirenberg inequality (\ref{1.13}), we obtain that $\int_{\R^2}|u_n|^4dx$ is also bounded uniformly in $n$. Furthermore, we deduce from (\ref{2.4}) that for sufficiently large $n>0$,
\begin{equation}\label{2.7}
\begin{split}
e(a)+1\geq& E_{a}(u_{n})\\
\geq&\frac{a^*-a}{2a^*}\int_{\R^2}\big|\nabla|u_{n}|\big|^2dx
-\frac{Ca^{3}a^*}{a^*-a}
+\int_{\R^2}\Big(V(x)-\frac{\Omega^2}{4}|x|^2\Big)|u_{n}|^2dx.
\end{split}
\end{equation}
In view of the above facts, we obtain that for any given $0<a<a^*$ and $0\leq\Omega<\Omega^*$, there exists a constant
$0<K=K(a,\Omega)<\infty$ independent of $n$, such that
\begin{equation}\label{2.8}
\begin{split}
\int_{\R^2}|u_n|^4dx<2K,\ \ \int_{\R^2}\Big(V(x)-\frac{\Omega^2}{4}|x|^2\Big)|u_n|^2dx<K \\
\int_{\R^2}\big|\nabla|u_n|\big|^2dx<K^2 \ \  \text{uniformly for large}\ \  n>0.
\end{split}
\end{equation}
We then deduce from (\ref{22}), (\ref{2.6}) and (\ref{2.8}) that for large $n>0$
\begin{equation}\label{2.9}
\begin{split}
&e(a)+\Big(1+\frac{C}{2}a^\frac{3}{2}\Big)K+1\\
\geq&e(a)+\frac{1}{2}\int_{\R^2}|u_{n}|^4dx+\frac{1}{2}\int_{\R^2}\int_{\R^2}
\ln\Big(1+\frac{1}{|x-y|}\Big)|u_{n}(x)|^2|u_{n}(y)|^2dxdy+1\\
\geq&e(a)+\frac{1}{2}\int_{\R^2}|u_{n}|^4dx-\frac{1}{2}\int_{\R^2}\int_{\R^2}
\ln|x-y||u_{n}(x)|^2|u_{n}(y)|^2dxdy+1\\
 \geq&\int_{\R^2}|\nabla u_n |^2dx-\Omega\int_{\R^2}x^{\bot}\cdot(iu_n, \nabla u_n)dx.
\end{split}
\end{equation}

Since $0\leq\Omega<\Omega^*$ is fixed, it follows from (\ref{2.3}) that
\begin{equation}\label{2.10}
\begin{split}
|x|^2\leq C(\Omega)\Big(V(x)-\frac{\Omega^2}{4}|x|^2\Big)\ \ \text{for sufficiently large} \ \ |x|>0.
\end{split}
\end{equation}
By the Cauchy-Schwarz inequality, we then obtain from (\ref{2.8}) and (\ref{2.10}) that for sufficiently large $R>0$ and $n>0$, there exists a constant $C(a,\Omega,R)>0$ such that
\begin{equation}\label{2.11}
\begin{aligned}
&\Big|\Omega\int_{\R^2}x^{\perp}\cdot(iu_{n},\nabla u_{n})dx\Big|\\
\leq&\frac{1}{2}\int_{\R^2}|\nabla u_n|^{2}dx+\frac{\Omega^2}{2}\int_{\R^2}|x|^2|u_n|^2dx\\
\leq&\frac{1}{2}\int_{\R^2}|\nabla u_n|^{2}dx+\frac{\Omega^2}{2}\Big
[C(R)\int_{B_{R}(0)}|u_n|^2dx\\
&+C(\Omega)\int_{B^{c}_{R}(0)}
\Big(V(x)-\frac{\Omega^{2}}{4}|x|^{2}\Big)|u_n|^2dx\Big]\\
\leq&\frac{1}{2}\int_{\R^2}|\nabla u_n|^{2}dx+C(a,\Omega,R).
\end{aligned}
\end{equation}
Thus, we deduce from (\ref{2.9}) and (\ref{2.11}) that for sufficiently large $R>0$ and $n>0$,
$$e(a)+1+\Big(1+\frac{C}{2}a^\frac{3}{2}\Big)K\geq \frac{1}{2}\int_{\R^2}|\nabla u_n|^2dx-C(a,\Omega,R),$$
which further implies that for any given $0<a<a^*$ and $0\leq\Omega<\Omega^*$, $\int_{\R^2}|\nabla u_n|^2dx$ is bounded uniformly for large $n>0$. Therefore, $\{u_n\}$ is also bounded uniformly in $\mathcal{H}$. By the compactness Lemma \ref{lem:2.1}, we then obtain that, passing to a subsequence if necessary, there exists $u\in\mathcal{H}$ such that
$$u_n\rightharpoonup u\ \ \text{weakly in}\ \ \mathcal{H}\ \ \text{as}\ \ n\rightarrow \infty,$$
and
$$\ \ u_n\rightarrow u\ \ \text{strongly in}\ \ L^q(\R^2,\C)~(2\leq q<\infty)\ \ \text{as}\ \ n\rightarrow \infty,$$
this indicates that $\|u\|_{2}^{2}=a$. For the logarithmic convolution term, we note from (5) of Lemma \ref{lem:2.2} that
\begin{equation*}
\int_{\R^2}\int_{\R^2}{\rm ln}|x-y||u(x)|^2|u(y)|^2dxdy\leqslant\liminf_{n\to\infty}\int_{\R^2}\int_{\R^2}{\rm ln}|x-y||u_n(x)|^2|u_n(y)|^2dxdy.
\end{equation*}
By the weak lower semicontinuity, we conclude from above that
$$e(a)\leq E_a(u)\leq\liminf_{n\rightarrow\infty}E_a(u_{n})=e(a).$$
Thus, $E_a(u)=e(a)$. This implies that for any given $0\leq\Omega<\Omega^*$ and $0<a<a^*$, there exists at least one minimizer for $e(a)$, and the proof of Proposition \ref{pro2.4} is therefore complete.
\qed
\vskip 0.1truein

We next consider the following nonexistence result, combining with Proposition \ref{pro2.4}, which then completes the proof of Theorem \ref{thm:1.1}.

\begin{prop}\label{pro2.5}
Let $Q=Q(|x|)$ be the unique positive solution of (\ref{1.6}). Suppose $V(x)$ satisfies (\ref{1.10}) such that $\Omega^*>0$ exists, where $\Omega^*$ is defined as in (\ref{1.11}). Then there is no minimizer for $e(a)$, provided that either $0\leq\Omega<\Omega^*$ and $a\geq a^*$ or  $\Omega>\Omega^*$ and $a>0$.
Moreover, we have $\lim\limits_{a\nearrow a^*}e(a)=-\infty$ for any $\Omega\geq 0$.
\end{prop}
\noindent\textbf{Proof.}
Choose a nonnegative function $\varphi(x)\in C_{0}^{\infty}(\R^2)$ such that $\varphi(x)=1$ for $|x|\leq1$ and $\varphi(x)=0$ for $|x|\geq 2$. For any $\tau>0$, define
\begin{equation}\label{2.12}
\begin{aligned}
w_{\tau}(x)=C_\tau\frac{\tau a^{\frac{1}{2}}}{\|Q\|_{2}}\varphi(x-x_\tau)Q\big(\tau(x-x_\tau)\big)e^{i\Omega S(x)},
\end{aligned}
\end{equation}
where $x_\tau$ is to be determined later, $S(x)=\frac{1}{2}x\cdot x_{\tau}^{\perp}$ is chosen such that $\nabla S(x)=\frac{1}{2}x_{\tau}^{\perp}$ in $\R^2$, and $C_\tau>0$ is chosen such that $w_{\tau}\in S_{a}$, then $C_\tau>0$ satisfies
\begin{equation*}
\frac{1}{C_\tau^{2}}=\frac{1}{\|Q\|_{2}^{2}}\int_{\R^2}Q^2(x)\varphi^2\Big(\frac{x}{\tau}\Big)dx
=1+O(\tau^{-\infty})\ \ \text{as}\ \ \tau\rightarrow\infty.
\end{equation*}
Here and below we use $f(t)=O(t^{-\infty})$ to denote a function $f$ satisfying $\lim_{t\rightarrow\infty}|f(t)|t^s=0$ for all $s>0$.

Direct calculations yield that
\begin{equation}\label{2.13}
\begin{aligned}
E_{a}(w_{\tau})=&\int_{\R^2}(|\nabla w_{\tau}|^2+V(x)|w_{\tau}|^2)dx-\Omega\int_{\R^2}x^{\bot}\cdot(iw_{\tau}, \nabla w_{\tau})dx\\
&-\frac{1}{2}\int_{\R^2}|w_{\tau}(x)|^4dx+\frac{1}{2}\int_{\R^2}\int_{\R^2}\ln|x-y||w_{\tau}(x)|^2
|w_{\tau}(y)|^2dxdy\\
=&\int_{\R^2}|(\nabla -i\mathcal{A})w_{\tau}|^2dx-\frac{1}{2}\int_{\R^2}|w_{\tau}(x)|^4dx+
\int_{\R^2}V_{\Omega}(x)|w_{\tau}|^2dx\\
&+\frac{1}{2}\int_{\R^2}\int_{\R^2}\ln|x-y||w_{\tau}(x)|^2
|w_{\tau}(y)|^2dxdy\\
=&\frac{a\tau^{2}}{a^*}\Big\{\int_{\mathbb{R}^{2}}\Big|\nabla Q-\frac{i\Omega}{2\tau^{2}}Q(x)x^{\bot}\Big|^{2}dx-\frac{a}{2a^*}
\int_{\mathbb{R}^{2}}Q^{4}dx+O(\tau^{-\infty})\Big\}\\
&+\int_{\R^2}V_{\Omega}(x)|w_{\tau}|^2dx
+\frac{1}{2}\int_{\R^2}\int_{\R^2}\ln|x-y||w_{\tau}(x)|^2
|w_{\tau}(y)|^2dxdy\\
\leq&a\tau^{2}(1-\frac{a}{a^{*}})
+\frac{1}{2}\Big(\frac{a}{a^*}\Big)^{2}\int_{\R^2}\int_{\R^2}\ln|x-y|Q^2(x)Q^2(y)dxdy\\
&-\frac{a^{2}}{2}\ln \tau+
\frac{aC^{2}_{\tau}}{a^*}\int_{\mathbb{R}^{2}}V_{\Omega}\Big(\frac{x}{\tau}+x_\tau\Big)
\varphi^{2}\Big(\frac{x}{\tau}\Big)Q^{2}(x)dx+O(\tau^{-2})
\end{aligned}
\end{equation}
as $\tau\rightarrow\infty$. We are now ready to prove the nonexistence of minimizers as follows.
\vskip 0.1truein
(i). If $0\leq\Omega<\Omega^*$ and $a\geq a^*$, take $x_\tau=x_{0}$, where $x_{0}$ is a certain fixed point in $\mathbb{R}^{2}$. It then follows from \cite[Theorem 1.8]{Lieb1} that
\begin{equation}\label{a5}
\lim_{\tau\rightarrow\infty}\Big[\frac{aC_{\tau}^{2}}{a^{*}}
\int_{\R^2}V_{\Omega}\Big(\frac{x}{\tau}+x_\tau\Big)
\varphi^{2}\Big(\frac{x}{\tau}\Big)Q^{2}(x)dx\Big]=aV_{\Omega}(x_{0}).
\end{equation}
We finally conclude from (\ref{22}), (\ref{2.13}) and (\ref{a5}) that
$$e(a)\leq E_{a}(w_{\tau})\leq -\frac{a^{2}}{4}\ln \tau\rightarrow-\infty\ \ \text{as}\ \ \tau\rightarrow\infty.$$
This implies that $e(a)$ is unbounded from below, and the nonexistence of minimizers for $e(a)$ is therefore proved.
\vskip 0.1truein
(ii). If $\Omega>\Omega^*$ and $a>0$, the definition of $\Omega^*$ implies that $\inf_{x\in\R^2}V_{\Omega}(x)=-\infty$ and there exists a certain constant $M_0>0$ such that $V_{\Omega}(x)<M_0$ holds in $\mathbb{R}^{2}$. Moreover, similar to (\ref{2.10}), we can obtain that
\begin{equation}\label{2.16}
\begin{split}
V_{\Omega}(x)=V(x)-\frac{\Omega^2}{4}|x|^2\leq-C(\Omega)|x|^2
\ \ \text{for sufficiently large} \ \ |x|>0.
\end{split}
\end{equation}
For any fixed and sufficiently large $\tau>0$, take a point $x_\tau=\big(\tau\sqrt{2\tau}, \tau\sqrt{2\tau}\big)\in\R^2$, direct calculations yield that
\begin{equation}\label{2.15}
V_{\Omega}\Big(\frac{x}{\tau}+x_{\tau}\Big)
\leq-C(\Omega)\Big|\frac{x}{\tau}+x_{\tau}\Big|^2
\leq-2\tau^{2}
\end{equation}
for any $x\in B_{\tau}(0)$.

It then follows from (\ref{22}), (\ref{2.13}) and (\ref{2.15}) that
\begin{equation*}
\begin{aligned}
E_{a}(w_{\tau})
\leq&a\tau^{2}
+\frac{1}{2}\Big(\frac{a}{a^*}\Big)^{2}\int_{\R^2}\int_{\R^2}\ln\big(1+|x-y|\big)Q^2(x)Q^2(y)dxdy
-\frac{a^{2}}{2}\ln \tau\\
&+\frac{aC^{2}_{\tau}}{a^*}\int_{B_{2\tau}(0)}V_{\Omega}\Big(\frac{x}{\tau}+x_\tau\Big)
\varphi^{2}\Big(\frac{x}{\tau}\Big)Q^{2}(x)dx+O(\tau^{-2})\\
\leq& a\tau^{2}
+\frac{1}{2}\Big(\frac{a}{a^*}\Big)^{2}\int_{\R^2}\int_{\R^2}\ln\big(1+|x-y|\big)Q^2(x)Q^2(y)dxdy
-\frac{a^{2}}{2}\ln \tau\\
&+\frac{aC^{2}_{\tau}}{a^*}\int_{B_{\tau}(0)}V_{\Omega}\Big(\frac{x}{\tau}+x_\tau\Big)
Q^{2}(x)dx+CM_0+O(\tau^{-2})\\
\leq& \Big[a-\frac{2aC_{\tau}^{2}}{a^*}\int_{B_{\tau}(0)}Q^{2}(x)dx\Big]\tau^{2}
-\frac{a^{2}}{2}\ln \tau+C+O(\tau^{-2})\\
\leq&  -\frac{a\tau^2}{2}-\frac{a^{2}}{4}\ln \tau\rightarrow-\infty\ \ \text{as}\ \ \tau\rightarrow\infty.
\end{aligned}
\end{equation*}
The above estimate implies that $e(a)$ is also unbounded from below in this case. This therefore gives the nonexistence of minimizers for the case where $\Omega>\Omega^*$ and $a>0$.

Moreover, for $a\in(0, a^*)$ and $\Omega\geq0$, take $x_{\tau}=x_0$ and $\tau=(a^*-a)^{-\frac{1}{2}}$, it then follows from (\ref{2.13}) and (\ref{a5}) that
$$\lim_{a\nearrow a^*}e(a)\leq\lim_{a\nearrow a^*}E_{a}(w_{\tau})=-\infty.$$
This implies that $\lim_{a\nearrow a^*}e(a)=-\infty$ and the proof of Proposition \ref{pro2.5} is complete.
\qed

\section{Mass Concentration as $a\nearrow a^*$}

By employing blow-up analysis and energy estimates, the main purpose of this section is to establish Theorem \ref{thm:1.2} on the concentration behavior of minimizers for $e(a)$ as $a\nearrow a^*$. Inspired by \cite{GLL0, GLY}, we first analyze some properties of minimizers for $e(a)$ with $0<a<a^*$.
\begin{lem}\label{lem:3.1}
For any $0<a<a^*$, set
\begin{equation}\label{3.1}
\begin{aligned}
\varepsilon_{a}:=\Big(\int_{\mathbb{R}^{2}}\big|\nabla|u_a|\big|^2dx\Big)^{-\frac{1}{2}}>0,
\end{aligned}
\end{equation}
where $u_a$ is a minimizer of $e(a)$. Suppose that the assumptions of Theorem \ref{thm:1.2} hold. Then we have
\begin{enumerate}
\item $\varepsilon_a>0$ satisfies
\begin{equation}\label{00}
\varepsilon_a\to0\ \ \text{and}\ \ \mu_{a}\varepsilon_{a}^2\rightarrow-\frac{1}{a^*}\ \
\text{as}\ \ a\nearrow a^*.
\end{equation}

  \item Define
\begin{equation}\label{3.2}
w_a(x):=\varepsilon_a u_a(\varepsilon_a x+x_a)e^{-i(\frac{\varepsilon_a \Omega}{2}x\cdot x_a^\bot-\theta_{a})},
\end{equation}
where $x_a$ is a global maximal point of $|u_a|$ and $\theta_{a}\in[0,2\pi)$ is a proper constant. Then there exists a constant $\eta>0$, independent of $0<a<a^*$, such that
\begin{equation}\label{3.3}
\int_{B_2(0)}|w_a(x)|^2dx\geqslant\eta>0\ \ \text{as}\ \ a\nearrow a^*.
\end{equation}
\item $w_{a}(x)$ satisfies
\begin{equation}\label{3.4}
|w_a(x)|\to\frac{1}{\sqrt{a^*}} Q\Big(\frac{x}{\sqrt{a^*}}\Big)\ \ \text{in}\ \ H^1(\R^2)\ \ \text{as}\ \ a\nearrow a^*,
\end{equation}
\end{enumerate}
where $Q$ is the unique solution of (\ref{1.6}).
\end{lem}
\vskip 0.1truein
\noindent\textbf{Proof.} 
1. We first prove that $\varepsilon_{a}\rightarrow 0$ as $a\nearrow a^*$. Following (\ref{1.13}), (\ref{1.12}) and (\ref{2.6}), we derive that for $a\in(0,a^*)$,
\begin{equation}\label{3.5}
\begin{aligned}
e(a)=&E_a(u_a)\\
=&\int_{\mathbb{R}^2}|(\nabla-i\mathcal{A})u_a|^2dx+\frac{1}{2}\int_{\R^2}\int_{\R^2}\ln|x-y|
|u_a(x)|^2|u_a(y)|^2dxdy\\
&-\frac{1}{2}\int_{\R^2}|u_a(x)|^4dx+\int_{\mathbb{R}^2}V_{\Omega}(x)|u_a(x)|^2dx\\
\geq&\frac{a^*-a}{a^*}\int_{\mathbb{R}^2}\big|\nabla|u_a|\big|^2dx-\frac{1}{2}\int_{\R^2}\int_{\R^2}
\ln\Big(1+\frac{1}{|x-y|}\Big)|u_a(x)|^2|u_a(y)|^2dxdy\\
\geq&\frac{a^*-a}{a^*}\varepsilon_{a}^{-2}-\frac{Ca^{\frac{3}{2}}}{2}\varepsilon_{a}^{-1}
\geq-\frac{Ca^{\frac{3}{2}}}{2}\varepsilon_{a}^{-1}.
\end{aligned}
\end{equation}
This together with the fact $\lim_{a\nearrow a^*}e(a)=-\infty$ yield that $\varepsilon_{a}\rightarrow 0$ as $a\nearrow a^*$.

We then claim that $\mu_{a}\varepsilon_{a}^2\rightarrow -\frac{1}{a^*}$ as $a\nearrow a^*$. In fact, by the definition of $\varepsilon_{a}$, we have
\begin{equation}\label{3.6}
\begin{aligned}
\varepsilon_a^2 e(a)=& \varepsilon_a^2\Big[\int_{\mathbb{R}^2}|(\nabla-i\mathcal{A})u_a|^2dx+\frac{1}{2}\int_{\R^2}
\int_{\R^2}\ln|x-y||u_a(x)|^2|u_a(y)|^2dxdy\\
&-\frac{1}{2}\int_{\R^2}|u_a(x)|^4dx+\int_{\mathbb{R}^2}V_{\Omega}(x)|u_a(x)|^2dx\Big]\\
\geq&1+\varepsilon_a^2\int_{\mathbb{R}^2}V_{\Omega}(x)|u_a(x)|^2dx
-\frac{\varepsilon_a^2}{2}\int_{\R^2}|u_a(x)|^4dx\\
&+\frac{\varepsilon_a^2}{2}\int_{\R^2}\int_{\R^2}\ln\big(1+|x-y|\big)|u_a(x)|^2|u_a(y)|^2dxdy\\
&-\frac{\varepsilon_a^2}{2}\int_{\R^2}\int_{\R^2}\ln\Big(1+\frac{1}{|x-y|}\Big)|u_a(x)|^2|u_a(y)|^2dxdy
.
\end{aligned}
\end{equation}
Since $\varepsilon_{a}\rightarrow 0$ as $a\nearrow a^*$, we derive from (\ref{2.6}) that
\begin{equation}\label{3.7}
\begin{aligned}
0\leq\frac{\varepsilon_{a}^2}{2}\int_{\R^2}\int_{\R^2}\ln\Big(1+\frac{1}{|x-y|}\Big)|u_a(x)|^2|u_a(y)|^2dxdy
\leq Ca^{\frac{3}{2}}\varepsilon_{a}\rightarrow 0\ \ \text{as}\ \ a\nearrow a^*.
\end{aligned}
\end{equation}
Combining with the fact that $\lim_{a\nearrow a^*}e(a)=-\infty$, we deduce from (\ref{1.13}), (\ref{3.6}) and (\ref{3.7}) that
\begin{equation*}
\begin{aligned}
0\geq&\limsup_{a\nearrow a^*}\varepsilon_{a}^2e(a)\\
\geq&\limsup_{a\nearrow a^*}\Big[\Big(1-\frac{\varepsilon_{a}^2}{2}\int_{\R^2}|u_a(x)|^4dx\Big)
+\varepsilon_{a}^2\int_{\mathbb{R}^2}V_{\Omega}(x)|u_a(x)|^2dx\\
&+\frac{\varepsilon_a^2}{2}\int_{\R^2}\int_{\R^2}\ln\big(1+|x-y|\big)|u_a(x)|^2|u_a(y)|^2dxdy\Big]\\
&+\liminf_{a\nearrow a^*}\Big[-\frac{\varepsilon_a^2}{2}\int_{\R^2}\int_{\R^2}\ln\Big(1+\frac{1}{|x-y|}\Big)|u_a(x)|^2|u_a(y)|^2
dxdy\Big]\\
=&\limsup_{a\nearrow a^*}\Big[\Big(1-\frac{\varepsilon_{a}^2}{2}\int_{\R^2}|u_a(x)|^4dx\Big)
+\varepsilon_{a}^2\int_{\mathbb{R}^2}V_{\Omega}(x)|u_a(x)|^2dx\\
&+\frac{\varepsilon_a^2}{2}\int_{\R^2}\int_{\R^2}\ln\big(1+|x-y|\big)|u_a(x)|^2|u_a(y)|^2dxdy\Big]\\
\geq&\limsup_{a\nearrow a^*}\Big[\Big(1-\frac{a}{a^*}\Big)
+\varepsilon_{a}^2\int_{\mathbb{R}^2}V_{\Omega}(x)|u_a(x)|^2dx\\
&+\frac{\varepsilon_a^2}{2}\int_{\R^2}\int_{\R^2}\ln\big(1+|x-y|\big)|u_a(x)|^2|u_a(y)|^2dxdy\Big]\geq0,
\end{aligned}
\end{equation*}
which implies that
\begin{equation}\label{3.8}
\begin{aligned}
&\lim\limits_{a\nearrow a^*}\varepsilon_a^2\int_{\R^2}|u_a|^4dx=2,
\ \ \lim\limits_{a\nearrow a^*}\varepsilon_a^2\int_{\R^2}V_{\Omega}(x)|u_{a}|^2dx=0,\\
&\lim\limits_{a\nearrow a^*}\varepsilon_a^2\int_{\R^2}\int_{\R^2}\ln\big(1+|x-y|\big)|u_a(x)|^2|u_a(y)|^2dxdy=0,\ \ \lim\limits_{a\nearrow a^*}\varepsilon_{a}^2e(a)=0.
\end{aligned}
\end{equation}
Note that $u_a$ is a minimizer of $e(a)$ and it satisfies the Euler-Lagrange equation (\ref{1.18}), one can check that the Lagrange multiplier $\mu_{a}\in\mathbb{R}$ satisfies
\begin{equation}\label{3.9}
\begin{aligned}
a\mu_{a}=e(a)+\frac{1}{2}\int_{\R^2}\int_{\R^2}\ln|x-y||u_a(x)|^2|u_a(y)|^2dxdy
-\frac{1}{2}\int_{\R^2}|u_a|^4dx.
\end{aligned}
\end{equation}
We thus obtain from \eqref{3.7}--\eqref{3.9} that
\begin{equation}\label{3.10}
\mu_a\varepsilon_a^2\to-\frac1{a^*}\ \  {\rm as}\ \ a\nearrow a^*.
\end{equation}
\vskip 0.1truein
2. Denote $\bar{w}_{a}(x):=\varepsilon_a u_a(\varepsilon_a x+x_a)e^{-i(\frac{\varepsilon_a \Omega}{2}x\cdot x_a^\bot)}$ and $w_{a}(x):=\bar{w}_{a}(x)e^{i\theta_{a}}$, where the parameter $\theta_{a}\in[0, 2\pi)$ is chosen such that
\begin{equation}\label{3.11}
\Big\|w_{a}-\frac{1}{\sqrt{a^*}}Q\Big(\frac{x}{\sqrt{a^*}}\Big)\Big\|_{L^{2}(\mathbb{R}^2)}
=\min_{\theta\in[0, 2\pi)}\Big\|e^{i\theta}\bar{w}_{a}(x)-\frac{1}{\sqrt{a^*}}Q\Big(\frac{x}{\sqrt{a^*}}\Big)\Big\|
_{L^{2}(\mathbb{R}^2)}.
\end{equation}
Rewrite $w_{a}(x):=R_{a}(x)+iI_{a}(x)$, where $R_{a}(x)$ and $I_{a}(x)$ denote the real and imaginary parts of $w_{a}(x)$, respectively. From (\ref{3.11}), we obtain the following orthogonality condition:
\begin{equation}\label{3.12}
\int_{\mathbb{R}^2}Q\Big(\frac{x}{\sqrt{a^*}}\Big)I_{a}(x)dx=0.
\end{equation}
Following (\ref{1.18}), we see that $w_{a}(x)$ satisfies the following Euler-Lagrange equation:
\begin{equation}\label{3.13}
\begin{aligned}
-\Delta w_{a}(x)+i\varepsilon_{a}^2\Omega\big(x^{\bot}\cdot\nabla w_{a}(x)\big)+\Big[\frac{\varepsilon_{a}^4\Omega^{2}}{4}|x|^{2}+\varepsilon_{a}^2V_{\Omega}
(\varepsilon_{a}x+x_{a})+a\varepsilon_{a}^{2}\ln\varepsilon_{a}\\
+\varepsilon_{a}^{2}\int_{\R^2}\ln|x-y||w_a(y)|^2dy
-\mu_{a}\varepsilon_{a}^{2}-|w_{a}(x)|^2\Big]w_{a}(x)=0
\ \ \text{in}\ \ \mathbb{R}^2.
\end{aligned}
\end{equation}
Define $W_{a}(x)=|w_{a}(x)|^2\geq0$ in $\mathbb{R}^2$, we then derive from (\ref{3.13}) that
\begin{equation}\label{3.14}
\begin{aligned}
-\frac{1}{2}\Delta W_{a}(x)+|\nabla w_{a}|^2-\varepsilon_{a}^{2}\Omega x^\bot\cdot(iw_{a}, \nabla w_{a})
+\Big[\frac{\varepsilon_{a}^4\Omega^{2}}{4}|x|^2+\varepsilon_{a}^2V_{\Omega}
(\varepsilon_{a}x+x_{a})\\+a\varepsilon_{a}^{2}\ln\varepsilon_{a}
+\varepsilon_{a}^2\int_{\R^2}\ln|x-y|W_{a}(y)dy-\mu_{a}\varepsilon_{a}^{2}-W_{a}\Big]W_{a}=0
\ \ \text{in}\ \ \mathbb{R}^2.
\end{aligned}
\end{equation}
By the diamagnetic inequality (\ref{1.12}), we have
$$|\nabla w_{a}|^2-\varepsilon_{a}^2\Omega x^{\bot}\cdot(iw_{a}, \nabla w_{a})+\frac{\varepsilon_{a}^4\Omega^{2}}{4}|x|^2|w_a|^2\geq0\ \ \text{in}\ \ \mathbb{R}^2,$$
together with (\ref{3.14}), which yields that
\begin{equation}\label{3.15}
\begin{aligned}
-\frac{1}{2}\Delta W_{a}(x)&+
\varepsilon_{a}^2\int_{\R^2}\ln|x-y|W_{a}(y)dyW_{a}+a\varepsilon_{a}^{2}\ln\varepsilon_{a}W_{a}\\
&-\mu_{a}\varepsilon_{a}^{2}W_{a}-W_{a}^{2}\leq0
\ \ \text{in}\ \ \mathbb{R}^2.
\end{aligned}
\end{equation}

On the other hand, since
\begin{equation}\label{eq216}
\int_{\R^2}|w_{a}|^{2}dx=a \ \ \text{and}\ \ \int_{\R^2}\big|\nabla |w_a|\big|^{2}dx=1,
\end{equation}
which imply that $\{|w_a|\}$ is bounded uniformly in $H^1(\mathbb{R}^2)$ for $a\in(0,a^*)$.
Direct calculations show that there exists a constant $C>0$ such that for any $x\in\mathbb{R}^{2}$ and $a\in(0,a^*)$
\begin{equation}\label{3.16}
\begin{aligned}
&\int_{\mathbb{R}^2}\ln\Big(1+\frac{1}{|x-y|}\Big)W_{a}(y)dy\\
=&\int_{\mathbb{R}^2}\ln\Big(1+\frac{1}{|x-y|}\Big)|w_{a}(y)|^2dy
\leq\int_{\mathbb{R}^2}\frac{1}{|x-y|}|w_{a}(y)|^2dy\\
=&\int_{|x-y|<1}\frac{1}{|x-y|}|w_{a}(y)|^2dy+\int_{|x-y|\geq1}\frac{1}{|x-y|}|w_{a}(y)|^2dy\\
\leq&\Big(\int_{|x-y|<1}\frac{1}{|x-y|^\frac{3}{2}}dy\Big)^{\frac{2}{3}}
\Big(\int_{|x-y|<1}|w_{a}(y)|^6dy\Big)^{\frac{1}{3}}+\int_{|x-y|\geq1}|w_{a}(y)|^2dy\\
\leq& C\big\||w_{a}|\big\|_{H^1(\mathbb{R}^2)}^{2}\leq C.
\end{aligned}
\end{equation}
We then derive from (\ref{3.10}), (\ref{3.15}) and (\ref{3.16}) that
\begin{equation}\label{3.17}
\begin{aligned}
-\Delta W_{a}+\frac{1}{2a^*}W_{a}-2W_{a}^{2}\leq0\ \ \text{in}\ \ \mathbb{R}^2\ \ \text{as}\ \ a\nearrow a^*.
\end{aligned}
\end{equation}
Since $0$ is a global maximal point of $W_{a}$ for all $0<a<a^*$, this implies that $-\Delta W_{a}(0)\geq0$ for $a\in(0, a^*)$. We then infer from (\ref{3.17}) that there exists some constant $\alpha>0$, independent of $a$, such that $W_{a}(0)\geq\alpha>0$ as $a\nearrow a^*$. Applying the De Giorgi-Nash-Moser theory \cite[Theorem 4.1]{Han}, we deduce from \eqref{3.17} that there exists a constant $C>0$ such that
\begin{equation}\label{3.18}
\begin{aligned}
\int_{B_{2}(0)}W_{a}^2dx\geq C\max_{x\in B_{1}(0)}W_{a}(x)\geq C_{1}(\alpha)\ \ \text{as}\ \ a\nearrow a^*.
\end{aligned}
\end{equation}
We then derive from (\ref{eq216}) and (\ref{3.18}) that
\begin{equation}\label{3.19}
\begin{aligned}
W_{a}(0)=\max_{x\in \mathbb{R}^{2}}W_{a}(x)\leq C\int_{B_{2}(0)}W_{a}^2dx\leq C\ \ \text{as}\ \ a\nearrow a^*.
\end{aligned}
\end{equation}
We finally conclude from (\ref{3.18}) and (\ref{3.19}) that
\begin{equation*}
\begin{aligned}
\int_{B_{2}(0)}|w_{a}|^2dx\geq\int_{B_{2}(0)}\frac{W_{a}^2(x)}
{\max_{x\in \mathbb{R}^{2}}W_{a}(x)}dx\geq C_{2}(\alpha):=\eta>0\ \ \text{as}\ \ a\nearrow a^*,
\end{aligned}
\end{equation*}
which completes the proof of (\ref{3.3}).
\vskip 0.1truein
3. Since $\{|w_a|\}$ is bounded uniformly in $H^1(\mathbb{R}^2)$, there exists a subsequence, still denoted by $\{|w_a|\}$, of $\{|w_a|\}$ such that $|w_a|\rightharpoonup w_0$ weakly in $H^1(\mathbb{R}^2)$ as $a\nearrow a^*$ for some $w_0\in H^1(\mathbb{R}^2)$, note from (\ref{3.3}) that $w_{0}\not\equiv 0$ in $\mathbb{R}^2$. By the weak convergence, we may assume that $|w_a|\rightarrow w_{0}$ almost everywhere in $\mathbb{R}^2$ as $a\nearrow a^*$. Applying the Br\'{e}zis-Lieb lemma \cite{Brezis}, we have
\begin{equation}\label{3.20}
\begin{aligned}
\|w_a\|_q^q=\|w_0\|_q^q+\big\||w_a|-w_0\big\|_q^q+o(1)\ \ {\rm as}\ \ a\nearrow a^*,\ \ {\rm where}\ \ 2\leq q<\infty,
\end{aligned}
\end{equation}
and
\begin{equation}\label{3.21}
\begin{aligned}
\big\|\nabla |w_a|\big\|_2^2=\|\nabla w_0\|_2^2+\big\|\nabla (|w_a|-w_0)\big\|_2^2+o(1)\ \ {\rm as}\ \ a\nearrow a^*,
\end{aligned}
\end{equation}
which imply that $\big\||w_a|-w_0\big\|_2^2< a^*$ as $a\nearrow a^*$.

Combining (\ref{1.13}) and (\ref{3.8}), we deduce from (\ref{3.20}) and (\ref{3.21}) that
\begin{equation}\label{3.22}
\begin{aligned}
0=&\lim_{a\nearrow a^*}\Big(\int_{\R^2}\big|\nabla |w_a|\big|^2dx-\frac12\int_{\R^2}|w_a|^4dx\Big)\\
=&\lim_{a\nearrow a^*}\Big(\int_{\R^2}\big|\nabla |w_a|-\nabla w_0\big|^2dx-\frac12\int_{\R^2}\big||w_a|-w_0\big|^4dx\Big)\\
&+\int_{\R^2}|\nabla w_0|^2dx-\frac12\int_{\R^2} |w_0|^4dx\\
\geqslant&\lim_{a\nearrow a^*}\Big(\int_{\R^2}\big|\nabla |w_a|-\nabla w_0\big|^2dx-\frac12\int_{\R^2}\big||w_a|-w_0\big|^4dx\Big)\\
&+\frac{1}{2}\Big(a^*\|w_{0}\|_{2}^{-2}-1\Big)\int_{\R^2}\left| w_0\right|^4dx\\
\geqslant&\lim_{a\nearrow a^*}\Big(1-\frac1{a^*}\int_{\R^2}\big||w_a|-w_0\big|^2dx\Big)\int_{\R^2}\big|\nabla |w_a|-\nabla w_0\big|^2dx\geqslant0,
\end{aligned}
\end{equation}
from which we obtain that $\|w_0\|_{2}^{2}=a^*$ and
$$
\lim_{a\nearrow a^*}\int_{\mathbb{R}^2}\big|\nabla|w_a|-\nabla w_0\big|^2dx=0.$$
Furthermore, we derive that
\begin{equation}\label{32}
\begin{aligned}
\int_{\R^2}|\nabla w_0|^2dx=1,\  \
|w_a|\rightarrow w_0 \ \ \text{strongly in}\ \ H^1(\mathbb{R}^2)\ \ \text{as}\ \ a\nearrow a^*.
\end{aligned}
\end{equation}
On the other hand, it follows from (\ref{3.22}) that
$$\int_{\R^2}|\nabla w_0|^2dx=\frac12\int_{\R^2} |w_0|^4dx,$$
together with the fact that $\|w_0\|_{2}^{2}=a^*$, we obtain that the identity of the inequality (\ref{1.13}) is achieved by $w_0$. Combining the above facts, we derive from (\ref{32}) that there exists a point $y_0\in \mathbb{R}^2$ such that, up to a subsequence if necessary,
$$|w_a(x)|\rightarrow w_0(x)=\frac{1}{\sqrt{a^*}}Q\Big(\frac{|x+y_0|}{\sqrt{a^*}}\Big)\ \ \text{strongly in}\ \ H^1(\mathbb{R}^2)\ \ \text{as}\ \ a\nearrow a^*.$$

Since the origin is a global maximum point of $|w_a|$ for $a\in(0,a^*)$, it must be a global maximum point of $w_0$, which implies that $y_0=0$, and hence
\begin{equation}\label{3.23}
\begin{aligned}
|w_a(x)|\rightarrow \frac{1}{\sqrt{a^*}}Q\Big(\frac{x}{\sqrt{a^*}}\Big)\ \ \text{strongly in}\ \ H^1(\mathbb{R}^2)\ \ \text{as}\ \ a\nearrow a^*.
\end{aligned}
\end{equation}
Moreover, since the convergence (\ref{3.23}) is independent of the choice of the subsequence, we deduce that (\ref{3.23}) holds essentially true for the whole sequence. This therefore completes the proof of Lemma \ref{lem:3.1}.
\qed
\vskip 0.1truein
We then establish the following crucial lemma, which is a weak version of Theorem \ref{thm:1.2}.

\begin{lem}\label{pro:3.2}
Under the assumptions of Theorem \ref{thm:1.2}, let $u_a$ be a minimizer of $e(a)$ for $a\in(0,a^*)$, and consider the sequences $\{w_a\}$ and $\{x_a\}$ defined in Lemma \ref{lem:3.1}. Then we have
\begin{enumerate}
\item
There exists a large $R>0$ such that
\begin{equation}\label{3.24}
|w_a(x)|\leq Ce^{-\frac{|x|}{3\sqrt{a^*}}}\ \ \text{for}\ \ |x|\geq R\ \ \text{as}\ \ a\nearrow a^*,
\end{equation}
where the constant $C>0$ is independent of $a\in(0,a^*)$.

\item
The function $w_a(x)$ satisfies
\begin{equation}\label{3.25}
w_a(x)\rightarrow\frac{1}{\sqrt{a^*}}Q\Big(\frac{x}{\sqrt{a^*}}\Big)\ \ \text{strongly in}\ \  H^1(\mathbb{R}^{2},\mathbb{C})\cap L^\infty(\mathbb{R}^{2},\mathbb{C})\ \ \text{as}\ \ a\nearrow a^*.
\end{equation}
Moreover, the global maximum point $x_a$ of $|u_a|$ is unique.
\end{enumerate}
\end{lem}
\vskip 0.1truein
\noindent\textbf{Proof.}  
1. Set $W_{a}=|w_{a}|^2$, then $W_{a}$ satisfies (\ref{3.14}). Recall from (\ref{3.23}) that $W_{a}(x)\rightarrow\frac{1}{a^*}Q^{2}(\frac{x}{\sqrt{a^*}})$ strongly in $L^2(\mathbb{R}^2)$ as $a\nearrow a^*$. Applying the De Giorgi-Nash-Moser theory \cite[Theorem 4.1]{Han} to (\ref{3.17}), we can obtain that there exists a sufficiently large $R>0$ such that
\begin{equation}\label{3.26}
W_{a}(x)\leq\frac{1}{36a^*}\ \ \text{uniformly for}\ \ |x|\geq R\ \ \text{as}\ \ a\nearrow a^*.
\end{equation}
We then derive from (\ref{3.17}) and (\ref{3.26}) that for sufficiently large $R>0$,
\begin{equation}\label{3.27}
\begin{aligned}
-\Delta W_{a}+\frac{4}{9a^*}W_{a}\leq0\ \ \text{uniformly for}\ \ |x|\geq R\ \ \text{as}\ \ a\nearrow a^*.
\end{aligned}
\end{equation}
Applying the comparison principle to (\ref{3.27}) gives that there exists a constant $C>0$ such that
\begin{equation}\label{3.28}
|w_a(x)|\leq Ce^{-\frac{|x|}{3\sqrt{a^*}}}\ \ \text{uniformly for}\ \ |x|\geq R\ \ \text{as}\ \ a\nearrow a^*.
\end{equation}
This completes the proof of \eqref{3.24}.
\vskip 0.1truein
2. We now prove that
\begin{equation}\label{3.35}
\begin{aligned}
w_a(x)\rightarrow\frac{1}{\sqrt{a^*}}Q\Big(\frac{x}{\sqrt{a^*}}\Big)\ \ \text{strongly in}\ \ H^1(\mathbb{R}^{2},\mathbb{C})\ \ \text{as}\ \ a\nearrow a^*.
\end{aligned}
\end{equation}
We first claim that $\{w_a\}$ is bounded uniformly in $H^1(\mathbb{R}^2,\mathbb{C})$ as $a\nearrow a^*$. In fact, by the definition of $w_a$ in (\ref{3.2}), we only need to prove that there exists a constant $C>0$, independent of $a$, such that
\begin{equation}\label{3.33}
\begin{aligned}
\int_{\mathbb{R}^2}|\nabla w_a|^2dx\leq C\  \ \text{as} \ \ a\nearrow a^*.
\end{aligned}
\end{equation}

Note that $\big\{|w_a|\big\}$ is bounded uniformly in $H^{1}(\mathbb{R}^2)$, multiplying (\ref{3.13}) by $\bar{w}_{a}(x)$ and integrating over $\mathbb{R}^2$, it follows from the Cauchy-Schwarz inequality and (\ref{2.6}) that
\begin{equation}\label{3.34}
\begin{aligned}
\int_{\mathbb{R}^2}|\nabla w_a|^2dx=&\varepsilon_{a}^{2}\Omega\int_{\mathbb{R}^{2}}x^{\perp}
\cdot(iw_{a},\nabla w_{a})dx-\frac{\varepsilon_{a}^{4}\Omega^{2}}{4}\int_{\mathbb{R}^{2}}
|x|^2|w_{a}|^2dx+\int_{\mathbb{R}^{2}}|w_{a}|^4dx\\
&-\varepsilon_{a}^{2}\int_{\mathbb{R}^{2}}\int_{\mathbb{R}^{2}}
\ln|x-y||w_{a}(x)|^2|w_{a}(y)|^2dxdy+\mu_{a}\varepsilon_{a}^{2}
\int_{\mathbb{R}^{2}}|w_{a}|^2dx\\
&-\varepsilon_{a}^{2}\int_{\mathbb{R}^{2}} V_{\Omega}(\varepsilon_ax+x_a)|w_{a}|^2dx-a\varepsilon_{a}^{2}\ln\varepsilon_{a}
\int_{\mathbb{R}^{2}}|w_{a}|^2dx\\
\leq&\varepsilon_{a}^{2}\Omega\int_{\mathbb{R}^{2}}x^{\perp}
\cdot(iw_{a},\nabla w_{a})dx+\mu_{a}\varepsilon_{a}^{2}
\int_{\mathbb{R}^{2}}|w_{a}|^2dx\\
&-a\varepsilon_{a}^{2}\ln\varepsilon_{a}
\int_{\mathbb{R}^{2}}|w_{a}|^2dx
+\int_{\mathbb{R}^{2}}|w_{a}|^4dx\\
&+\varepsilon_{a}^{2}\int_{\mathbb{R}^{2}}
\int_{\mathbb{R}^{2}}\ln\Big(1+\frac{1}{|x-y|}\Big)|w_{a}(x)|^2|w_{a}(y)|^2dxdy\\
\leq&\frac{\varepsilon_{a}^{2}}{2}\int_{\mathbb{R}^{2}}|\nabla w_{a}|^2dx+
\frac{\varepsilon_{a}^{2}\Omega^{2}}{2}\int_{\mathbb{R}^{2}}
|x|^2|w_{a}|^2dx+a\mu_{a}\varepsilon_{a}^{2}\\
&+Ca^{\frac{3}{2}}\varepsilon_{a}^{2}\Big(\int_{\mathbb{R}^{2}}|\nabla w_{a}|^2dx\Big)^{\frac{1}{2}}+C+o(\varepsilon_{a})\\
\leq&\big(\frac{\varepsilon_{a}^{2}}{2}+\varepsilon_{a}\big)\int_{\mathbb{R}^{2}}|\nabla w_{a}|^2dx+a\mu_{a}\varepsilon_{a}^{2}+C+o(\varepsilon_{a})\ \ \text{as}\ \ a\nearrow a^*,
\end{aligned}
\end{equation}
where $C>0$ is independent of $a$. We then derive from (\ref{00}) and (\ref{3.34}) that (\ref{3.33}) holds true.

Since $\{w_a\}$ is bounded uniformly in $H^1(\mathbb{R}^{2},\mathbb{C})$, we may assume that up to a subsequence if necessary,
\begin{equation}\label{3.36}
\begin{aligned}
w_a\rightharpoonup w_1\ \ \text{weakly in}\ \ H^1(\mathbb{R}^{2},\mathbb{C})\ \ \text{as}\ \ a\nearrow a^*,
\end{aligned}
\end{equation}
and
\begin{equation}\label{3.37}
\begin{aligned}
w_a\rightarrow w_1\ \ \text{strongly in}\ \ L_{loc}^q(\mathbb{R}^{2},\mathbb{C})~(2\leq q<\infty)\ \ \text{as}\ \ a\nearrow a^*
\end{aligned}
\end{equation}
for some $w_1\in H^1(\mathbb{R}^{2},\mathbb{C})$ and $w_1\not\equiv 0$. From (\ref{3.4}) and (\ref{3.37}), we derive that
\begin{equation}\label{3.38}
\begin{aligned}
\int_{\mathbb{R}^2}|w_1|^2dx&=\lim_{R\rightarrow\infty}\int_{B_{R}(0)}|w_1|^2dx
=\lim_{R\rightarrow\infty}\lim_{a\nearrow a^*}\int_{B_{R}(0)}|w_a|^2dx\\
&=\lim_{R\rightarrow\infty}
\int_{B_{R}(0)}\Big|\frac{1}{\sqrt{a^*}}Q\Big(\frac{x}{\sqrt{a^*}}\Big)\Big|^2dx=a^*,
\end{aligned}
\end{equation}
which further implies that
\begin{equation}\label{3.39}
\begin{aligned}
w_a\rightarrow w_1\ \ \text{strongly in}\ \ L^2(\mathbb{R}^2,\mathbb{C})\ \ \text{as}\ \ a\nearrow a^*.
\end{aligned}
\end{equation}
Applying the interpolation inequality, we derive from (\ref{3.36}) and (\ref{3.39}) that
\begin{equation}\label{3.40}
\begin{aligned}
w_a\rightarrow w_1\ \ \text{strongly in}\ \ L^q(\mathbb{R}^2,\mathbb{C})~(2\leq q<\infty)\ \ \text{as}\ \ a\nearrow a^*.
\end{aligned}
\end{equation}
In view of (\ref{00}), (\ref{3.8}) and (\ref{3.40}), we obtain that
\begin{equation}\label{3.41}
\begin{aligned}
\int_{\mathbb{R}^2}|w_1|^4dx=\lim_{a\nearrow a^*}\int_{\mathbb{R}^2}|w_a|^4dx=2,\ \
\lim_{a\nearrow a^*}\mu_{a}\varepsilon_{a}^{2}
\int_{\mathbb{R}^{2}}|w_{a}|^2dx=-1.
\end{aligned}
\end{equation}
Combining \cite[Theorem 6.17]{Lieb1} with the Gagliardo-Nirenberg inequality (\ref{1.13}), we deduce from (\ref{3.34}), (\ref{3.36}) and (\ref{3.41}) that
$$1=\lim_{a\nearrow a^*}\int_{\mathbb{R}^2}|\nabla w_a|^2dx
\geq\int_{\mathbb{R}^2}|\nabla w_1|^2dx
\geq\int_{\mathbb{R}^2}\big|\nabla |w_1|\big|^2dx
\geq\frac{1}{2}\int_{\mathbb{R}^2}|w_1|^4dx=1,$$
which implies that
\begin{equation}\label{3.42}
\begin{aligned}
\lim_{a\nearrow a^*}\int_{\mathbb{R}^2}|\nabla w_a|^2dx
=\int_{\mathbb{R}^2}|\nabla w_1|^2dx=1,
\end{aligned}
\end{equation}
and
\begin{equation}\label{3.43}
\begin{aligned}
\int_{\mathbb{R}^2}|\nabla w_1|^2dx=\int_{\mathbb{R}^2}\big|\nabla |w_1|\big|^2dx,\ \ i.e.\ \ |\nabla w_1|=\big|\nabla |w_1|\big|\ \ a.e.\ \ \text{in}\ \ \mathbb{R}^2.
\end{aligned}
\end{equation}
Combining (\ref{3.36}), (\ref{3.39}), (\ref{3.42}) and (\ref{3.43}), we obtain from Lemma \ref{lem:3.1} that
\begin{equation}\label{3.44}
\begin{aligned}
\lim_{a\nearrow a^*}w_{a}(x)=\frac{1}{\sqrt{a^*}}Q\Big(\frac{x}{\sqrt{a^*}}\Big)e^{i\sigma}
\ \ \text{strongly in}\ \ H^1(\mathbb{R}^{2},\mathbb{C})
\end{aligned}
\end{equation}
for some $\sigma\in\mathbb{R}$. We further have $\sigma=0$ in (\ref{3.44}) in view of (\ref{3.12}), which implies that (\ref{3.35}) is proved.

Next, we prove that
\begin{equation}\label{3441}
\begin{aligned}
w_a(x)\rightarrow \frac{1}{\sqrt{a^*}}Q\Big(\frac{x}{\sqrt{a^*}}\Big)\ \ \text{uniformly in}\ \  L^\infty(\mathbb{R}^2,\mathbb{C})\ \ \text{as}\ \ a\nearrow a^*.
\end{aligned}
\end{equation}
Indeed, by the exponential decay of $Q(x)$ and $w_a(x)$, we only need to show the $L^\infty$--uniform convergence of $w_a(x)$ on any compact domain of $\mathbb{R}^2$ as $a\nearrow a^*$. Since $w_a(x)$ satisfies (\ref{3.13}), denote
\begin{equation*}
\begin{aligned}
G_{a}(x):=&-i\varepsilon_{a}^2\Omega\big(x^{\bot}\cdot\nabla w_{a}(x)\big)-\Big[\frac{\varepsilon_{a}^4\Omega^{2}}{4}|x|^{2}+\varepsilon_{a}^2V_{\Omega}
(\varepsilon_{a}x+x_{a})-\mu_{a}\varepsilon_{a}^{2}\\
&+a\varepsilon_{a}^{2}\ln\varepsilon_{a}
+\varepsilon_{a}^{2}\int_{\R^2}\ln|x-y||w_a(y)|^2dy
-|w_{a}(x)|^2\Big]w_{a}(x),
\end{aligned}
\end{equation*}
so that
\begin{equation}\label{3.45}
\begin{aligned}
-\Delta w_{a}(x)=G_{a}(x)\ \ \text{in}\ \ \mathbb{R}^2.
\end{aligned}
\end{equation}
Note that $\{w_a\}$ is bounded uniformly in $H^1(\mathbb{R}^2,\mathbb{C})$ and satisfies the exponential decay (\ref{3.24}), we deduce from (\ref{3.16}) that
\begin{equation}\label{3.451}
\begin{aligned}
&\Big|\int_{\R^2}\ln|x-y||w_a(y)|^2dy\Big|\\
\leq&\int_{\R^2}\ln\big(1+|x-y|\big)|w_a(y)|^2dy
+\int_{\R^2}\ln\Big(1+\frac{1}{|x-y|}\Big)|w_a(y)|^2dy\\
\leq&\int_{\R^2}\Big[\ln\big(1+|x|\big)+\ln\big(1+|y|\big)\Big]|w_a(y)|^2dy
+C\big\|w_{a}\big\|_{H^1(\mathbb{R}^2)}^{2}\\
\leq&a\ln\big(1+|x|\big)+C\ \ \text{as}\ \ a\nearrow a^*,
\end{aligned}
\end{equation}
where the constant $C>0$ is independent of $a$.

Combining the above facts, we obtain that $G_{a}(x)$ is also bounded uniformly in $L_{loc}^{2}(\mathbb{R}^2,\mathbb{C})$ as $a\nearrow a^*$. For any large $R>0$, it follows from \cite[Theorem 8.8]{Gilbarg} that
\begin{equation}\label{0}
\|w_a\|_{H^2(B_R)}\leqslant C\Big(\|w_a\|_{H^1(B_{R+1})}+\|G_a\|_{L^2(B_{R+1})}\Big),
\end{equation}
where $C>0$ is independent of $a$ and $R$. Therefore, $\{w_a\}$ is also bounded uniformly in $H^2(B_R)$ as $a\nearrow a^*$. Since the embedding $H^2(B_R)\hookrightarrow L^\infty(B_R)$ is compact, we obtain that there exists a subsequence $\{w_{a_{k}}\}$  of $\{w_a\}$ such that
\begin{equation*}
\lim_{a_{k}\nearrow a^*}w_{a_{k}}(x)= w_0(x)\ \ \text{uniformly in}\ \ L^\infty(B_R,\mathbb{C}).
\end{equation*}
Note that $R > 0$ is arbitrary, we get from (\ref{3.35}) that
\begin{equation}\label{3.46}
\lim_{a_{k}\nearrow a^*}w_{a_{k}}(x)=\frac{1}{\sqrt{a^*}}Q\Big(\frac{x}{\sqrt{a^*}}\Big)\ \ \text{uniformly in}\ \ L_{loc}^\infty(\R^2,\mathbb{C}).
\end{equation}
Since the above convergence is independent of the subsequence that we choose, (\ref{3.46}) holds essentially true for the whole sequence. This further implies that \eqref{3441} holds true.

Finally, we prove the uniqueness of $x_a$ as $a\nearrow a^*$, where $x_a$ is a global maximal point of $|u_a|$. It follows from (\ref{0}) that $\big\{|\nabla w_a|\big\}$ is bounded uniformly in $L_{loc}^{q}(\mathbb{R}^{2})$ as $a\nearrow a^*$ for any $q\geq 2$. The $L^p$ estimate \cite[Theorem 9.11]{Gilbarg} applied to (\ref{3.45}) then yields that $\{w_a\}$ is bounded uniformly in $W^{2,q}_{loc}(\R^2)$ as $a\nearrow a^*$. The standard Sobolev embedding thus gives that $\{w_a\}$ is bounded uniformly in $C^{1, \alpha}_{loc}(\R^2)$ as $a\nearrow a^*$. On the other hand, Lemma \ref{lem:A.3} in the Appendix shows that
$\varepsilon_a^2\int_{\R^2}\ln|x-y||w_a(y)|^2dy\in C^{\alpha}_{loc}(\R^2)$ and
$\big\{\varepsilon_a^2\int_{\R^2}\ln|x-y||w_a(y)|^2dy\big\}$ is bounded uniformly in $C^{ \alpha}_{loc}(\R^2)$ as $a\nearrow a^*$.
Furthermore, since $\varepsilon_a^2V_{\Omega}(\varepsilon_a x+x_a)\in C_{loc}^{\alpha}(\mathbb{R}^2)$, we have $\{G_{a}\}$ is bounded uniformly in $C^{ \alpha}_{loc}(\R^2)$ as $a\nearrow a^*$. It then follows from the Schauder estimate \cite[Theorem 6.2]{Gilbarg} that $\{w_a\}$ is bounded uniformly in $C^{2,\alpha}_{loc}(\R^2)$ as $a\nearrow a^*$. Therefore, up to a subsequence if necessary, there exists $\tilde{w}_{0}\in C_{loc}^{2}(\mathbb{R}^2)$ such that
$$w_a\rightarrow \tilde{w}_0\ \ \text{in}\ \ C^{2}_{loc}(\R^2)\ \ \text{as}\ \ a\nearrow a^*,$$
and we further deduce from (\ref{3.25}) that
\begin{equation}\label{3.47}
|w_a(x)|\rightarrow \tilde{w}_0(x)=\frac{1}{\sqrt{a^*}}Q\Big(\frac{x}{\sqrt{a^*}}\Big)\ \ \text{in}\ \ C^{2}_{loc}(\R^2)\ \ \text{as}\ \ a\nearrow a^*.
\end{equation}

Since the origin is the unique global maximal point of $Q(x)$, \eqref{3.47} shows that all global maximal points of $|w_a(x)|$ must stay in a small ball $B_\delta(0)$ as $a\nearrow a^*$ for some $\delta>0$. Noting $Q''(0)<0$, we conclude that $Q''(r)<0$ for $0\leq r<\delta$. It then follows from \eqref{3.47} and \cite[Lemma 4.2]{Ni} that each $|w_a|$ has a unique maximum point as $a\nearrow a^*$, which is exactly the origin. This further proves the uniqueness of global maximal points for $|u_a|$ as $a\nearrow a^*$ and the proof of Lemma \ref{pro:3.2} is complete.
\qed
\vskip 0.1truein

Based on Lemma \ref{lem:3.1} and Lemma \ref{pro:3.2}, the rest part of this section is to complete the proof of Theorem \ref{thm:1.2}. In order to reach this aim, we shall first establish the refined energy estimates of $e(a)$ and then analyze the blow-up point $x_a$ of $u_a$.
\vskip 0.1truein
\noindent\textbf{Completion of the proof for Theorem \ref{thm:1.2}.}
In view of (\ref{3.2}) and (\ref{3.25}), to establish the proof of Theorem \ref{thm:1.2}, it remains to prove that
\begin{equation}\label{3.48}
\begin{aligned}
\varepsilon_a&:=\Big(\int_{\mathbb{R}^{2}}\big|\nabla|u_a|\big|^2dx\Big)^{-\frac{1}{2}}\\
&=\frac{2}{a^*}\Big(\frac{a^*-a}{a^*}\Big)^{\frac{1}{2}}+o\Big((a^*-a)^{\frac{1}{2}}\Big)>0\ \
 \hbox{as}\ \ a\nearrow a^*,
\end{aligned}
\end{equation}
and $\lim_{a\nearrow a^*}V_{\Omega}(x_a)=0$.

In order to prove (\ref{3.48}), we consider the  trial function
$$u_\tau(x)=\frac{\tau a^{\frac{1}{2}}}{\|Q\|_{2}}Q(\tau x),\ \ \tau>0.$$
By the definition of $u_\tau(x)$, some calculations yield that
\begin{equation}\label{3.53}
\begin{aligned}
\int_{\mathbb{R}^{2}}\Big|\Big(\nabla-\frac{i\Omega}{2}x^{\bot}\Big)u_\tau(x)\Big|^{2}dx
&=\frac{a\tau^{2}}{\|Q\|_{2}^{2}}
\int_{\mathbb{R}^{2}}\Big|\Big(\nabla-\frac{i\Omega}{2}x^{\bot}\Big)Q(\tau x)\Big|^{2}dx\\
&=\frac{a}{a^*}\int_{\mathbb{R}^{2}}\Big|\tau\nabla Q(x)-\frac{i\Omega x^{\bot}}{2\tau}Q(x)\Big|^{2}dx\\
&=a\tau^{2}+\frac{a\Omega^{2}}{4a^*\tau^{2}}\int_{\mathbb{R}^{2}}|x|^2Q^2(x)dx,
\end{aligned}
\end{equation}
and
\begin{equation}\label{3.54}
\begin{aligned}
\int_{\mathbb{R}^{2}}V_{\Omega}(x)|u_\tau(x)|^{2}dx
&=\frac{a\tau^{2}}{\|Q\|_{2}^{2}}\int_{\mathbb{R}^{2}}
V_{\Omega}(x)\big|Q(\tau x)\big|^{2}dx\\
&=\frac{a}{a^*}\int_{\mathbb{R}^{2}}
V_{\Omega}\Big(\frac{x}{\tau}\Big)Q^{2}(x)dx.
\end{aligned}
\end{equation}
We also deduce that
\begin{equation}\label{3.55}
\begin{aligned}
&\frac{1}{2}\int_{\mathbb{R}^{2}}\int_{\mathbb{R}^{2}}\ln|x-y||u_\tau(x)|^2|u_\tau(y)|^2dxdy\\
=&\frac{a^{2}\tau^{4}}{2(a^*)^{2}}\int_{\mathbb{R}^{2}}\int_{\mathbb{R}^{2}}
\ln|x-y|Q^{2}(\tau x)Q^{2}(\tau y)dxdy\\
=&\frac{1}{2}\Big(\frac{a}{a^*}\Big)^{2}\int_{\mathbb{R}^{2}}\int_{\mathbb{R}^{2}}
\ln|x-y|Q^{2}(x)Q^{2}(y)dxdy-\frac{a^{2}}{2}\ln\tau,
\end{aligned}
\end{equation}
and
\begin{equation}\label{3.56}
\begin{aligned}
\frac{1}{2}\int_{\mathbb{R}^{2}}|u_{\tau}(x)|^{4}dx
=\frac{a^{2}\tau^{4}}{2(a^*)^{2}}\int_{\mathbb{R}^{2}}Q^{4}(\tau x)dx
=\frac{a^{2}\tau^{2}}{a^*}.
\end{aligned}
\end{equation}
It then follows from (\ref{3.53})--(\ref{3.56}) that
\begin{equation}\label{3.57}
\begin{aligned}
e(a)\leq& E_{a}\big(u_{\tau}(x)\big)\\
=&\int_{\mathbb{R}^{2}}\Big|\Big(\nabla-\frac{i\Omega}{2}x^{\bot}\Big)u_\tau(x)\Big|^{2}dx
+\int_{\mathbb{R}^{2}}V_{\Omega}(x)|u_\tau(x)|^{2}dx
-\frac{1}{2}\int_{\mathbb{R}^{2}}|u_{\tau}(x)|^{4}dx\\
&+\frac{1}{2}\int_{\mathbb{R}^{2}}\int_{\mathbb{R}^{2}}\ln|x-y||u_\tau(x)|^2|u_\tau(y)|^2dxdy\\
=&\frac{a(a^{*}-a)}{a^*}\tau^{2}+\frac{a\Omega^{2}}{4a^*\tau^{2}}\int_{\mathbb{R}^{2}}|x|^2Q^2(x)dx
+\frac{a}{a^*}\int_{\mathbb{R}^{2}}
V_{\Omega}\Big(\frac{x}{\tau}\Big)Q^{2}(x)dx\\
&+\frac{1}{2}\Big(\frac{a}{a^*}\Big)^{2}\int_{\mathbb{R}^{2}}\int_{\mathbb{R}^{2}}
\ln|x-y|Q^{2}(x)Q^{2}(y)dxdy-\frac{a^{2}}{2}\ln\tau.
\end{aligned}
\end{equation}
Setting
$$\tau=\Big[\frac{aa^*}{4(a^*-a)}\Big]^{\frac{1}{2}}>0,$$
we deduce from ($V$) and (\ref{3.57}) that
\begin{equation}\label{3.58}
\begin{aligned}
e(a)\leq&\frac{(a^*)^{2}}{4}-\frac{(a^*)^{2}}{2}\ln a^*+\frac{1}{2}\int_{\mathbb{R}^{2}}\int_{\mathbb{R}^{2}}
\ln|x-y|Q^2(x)Q^2(y)dxdy\\
&+\frac{a^{2}}{4}\ln\big[4(a^*-a)\big]+o(1)\ \ \text{as}\ \ a\nearrow a^*.
\end{aligned}
\end{equation}

On the other hand, by the definition of $w_a(x)$ in (\ref{3.2}), we obtain that
\begin{equation}\label{3.49}
\begin{aligned}
&\int_{\mathbb{R}^2}\Big|\Big(\nabla-\frac{i\Omega}{2}x^{\bot}\Big)u_a(x)\Big|^{2}dx\\
=&\int_{\mathbb{R}^2}\Big|\Big(\nabla-\frac{i\Omega}{2}x^{\bot}\Big)\frac{1}{\varepsilon_{a}}
w_{a}\Big(\frac{x-x_{a}}{\varepsilon_{a}}\Big)
e^{i\big(\frac{\Omega}{2}x\cdot x_{a}^{\bot}-\theta_{a}\big)}\Big|^{2}dx\\
=&\int_{\mathbb{R}^2}\Big|\frac{1}{\varepsilon_{a}}\nabla w_{a}(x)-\frac{i\varepsilon_{a}\Omega}{2}x^{\bot}w_{a}(x)\Big|^{2}dx\\
=&\frac{1}{\varepsilon_{a}^{2}}\int_{\mathbb{R}^2}|\nabla w_{a}|^{2}dx
+\frac{\varepsilon_{a}^{2}\Omega^{2}}{4}\int_{\mathbb{R}^2}|x|^2|w_{a}(x)|^{2}dx
-\Omega\int_{\mathbb{R}^2}x^{\bot}\cdot(iw_{a},\nabla w_{a})dx.
\end{aligned}
\end{equation}
Applying the exponential decay (\ref{3.24}) and the convergence (\ref{3.25}), we have
\begin{equation}\label{3.50}
\begin{aligned}
&\Omega\Big|\int_{\mathbb{R}^{2}}x^{\bot}\cdot(iw_{a},\nabla w_{a})dx\Big|\\
=&\Omega\Big|\int_{\mathbb{R}^{2}}x^{\bot}\cdot(R_{a}\nabla I_{a}-I_{a}\nabla R_{a})dx\Big|\\
=&2\Omega\Big|\int_{\mathbb{R}^{2}}x^{\bot}\cdot(R_{a}\nabla I_{a})dx\Big|\leq C\|\nabla I_{a}\|_{L^2(\mathbb{R}^{2})}=o(1)\ \ \text{as}\ \ a\nearrow a^*.
\end{aligned}
\end{equation}
Furthermore, it follows from (\ref{3.25}), (\ref{3.41}), (\ref{3.49}) and (\ref{3.50}) that
\begin{equation}\label{3.51}
\begin{aligned}
e(a)=&E_{a}(u_a)=E_{a}\Big(\frac{1}{\varepsilon_{a}}
w_{a}\Big(\frac{x-x_a}{\varepsilon_{a}}\Big)e^{i\big(\frac{\Omega}{2}x\cdot x_{a}^{\bot}-\theta_{a}\big)}\Big)\\
=&\frac{1}{\varepsilon_{a}^{2}}\int_{\mathbb{R}^2}|\nabla w_{a}|^{2}dx
+\frac{\varepsilon_{a}^{2}\Omega^{2}}{4}\int_{\mathbb{R}^2}|x|^2|w_{a}|^{2}dx
-\Omega\int_{\mathbb{R}^2}x^{\bot}\cdot(iw_{a},\nabla w_{a})dx\\
&-\frac{1}{2\varepsilon_{a}^{2}}\int_{\mathbb{R}^2}|w_{a}|^4dx
+\int_{\mathbb{R}^2}V_{\Omega}(\varepsilon_{a}x+x_a)|w_{a}(x)|^{2}dx
+\frac{a^2}{2}\ln\varepsilon_{a}\\
&+\frac{1}{2}\int_{\mathbb{R}^{2}}\int_{\mathbb{R}^{2}}
\ln|x-y||w_{a}(x)|^2|w_{a}(y)|^2dxdy\\
\geq&\frac{1}{\varepsilon_{a}^{2}}\Big[\int_{\mathbb{R}^2}|\nabla w_a|^2dx-\frac{a^*}{2a}\int_{\mathbb{R}^2}|w_a|^4dx\Big]
+\frac{a^2}{2}\ln\varepsilon_{a}\\
&+\frac{1}{2\varepsilon_{a}^{2}}\big(\frac{a^*}{a}-1\big)\int_{\mathbb{R}^2}|w_a|^4dx
-\Omega\int_{\mathbb{R}^2}x^{\bot}\cdot(iw_{a},\nabla w_{a})dx\\
&+\frac{1}{2}\int_{\mathbb{R}^{2}}\int_{\mathbb{R}^{2}}
\ln|x-y||w_{a}(x)|^2|w_{a}(y)|^2dxdy\\
\geq&\frac{a^*-a}{2a\varepsilon_{a}^{2}}\int_{\mathbb{R}^2}|w_a|^4dx
+\frac{a^2}{2}\ln\varepsilon_{a}-\Omega\int_{\mathbb{R}^2}x^{\bot}\cdot(iw_{a},\nabla w_{a})dx\\
&+\frac{1}{2}\int_{\mathbb{R}^{2}}\int_{\mathbb{R}^{2}}\ln|x-y||w_{a}(x)|^2|w_{a}(y)|^2dxdy
\\
\geq&\big[1+o(1)\big]\frac{a^*-a}{a\varepsilon_{a}^{2}}
+\frac{1}{2}\int_{\mathbb{R}^{2}}\int_{\mathbb{R}^{2}}
\ln|x-y|Q^2(x)Q^2(y)dxdy\\
&+\frac{a^2}{2}\ln\varepsilon_{a}+\frac{(a^*)^{2}}{4}\ln a^*+o(1)\\
\geq&\frac{(a^*)^{2}}{4}-\frac{(a^*)^{2}}{2}\ln a^*+\frac{1}{2}\int_{\mathbb{R}^{2}}\int_{\mathbb{R}^{2}}
\ln|x-y|Q^2(x)Q^2(y)dxdy\\
&+\frac{a^{2}}{4}\ln\big[4(a^*-a)\big]+o(1)\ \ \text{as}\ \ a\nearrow a^*.
\end{aligned}
\end{equation}
Combining this with the upper energy estimate (\ref{3.58}), we obtain that
\begin{equation}\label{3.59}
\begin{aligned}
e(a)\approx&\frac{(a^*)^{2}}{4}-\frac{(a^*)^{2}}{2}\ln a^*+\frac{1}{2}\int_{\mathbb{R}^{2}}\int_{\mathbb{R}^{2}}
\ln|x-y|Q^2(x)Q^2(y)dxdy\\
&+\frac{a^{2}}{4}\ln\big[4(a^*-a)\big]\ \ \text{as}\ \ a\nearrow a^*,
\end{aligned}
\end{equation}
and the last equality of (\ref{3.51}) holds if and only if $\varepsilon_{a}>0$ satisfies (\ref{3.48}).

Finally, we prove that $\lim_{a\nearrow a^*}V_{\Omega}(x_a)=0$. Motivated by \cite{GZZ}, we claim that $\{x_a\}$ is bounded uniformly for $a\nearrow a^*$. On the contrary, we assume that there exists a subsequence, denoted
by $\{x_{a_{k}}\}$ with $a_{k}\nearrow a^*$ as $k\rightarrow \infty$, such that
$$x_{a_{k}}\rightarrow\infty \ \ \text{as} \ \ k\rightarrow \infty.$$
By the definition of $\Omega^*$, there exists $C_0>0$ such that
\begin{equation}\label{3.600}
\liminf_{k\rightarrow\infty} V_{\Omega}(\varepsilon_{a_{k}}x+x_{a_{k}})>C_0\ \ \text{for}\ \ x\in B_{2}(0).
\end{equation}
Moreover, we derive from (\ref{3.3}), (\ref{3.51}) (\ref{3.59}) and (\ref{3.600}) that
\begin{equation}\label{3.60}
\begin{aligned}
0&=\lim_{a\nearrow a^*}\int_{\mathbb{R}^2}V_{\Omega}(\varepsilon_{a}x+x_a)|w_{a}(x)|^{2}dx\\
&=\lim_{k\rightarrow\infty}\int_{\mathbb{R}^2}V_{\Omega}(\varepsilon_{a_{k}}x+x_{a_{k}})
|w_{a_{k}}(x)|^{2}dx\\
&\geq\liminf_{k\rightarrow\infty}\int_{B_{2}(0)}V_{\Omega}(\varepsilon_{a_{k}}x+x_{a_{k}})
|w_{a_{k}}(x)|^{2}dx\\
&\geq \frac{C_{0}\eta}{2}>0,
\end{aligned}
\end{equation}
which is impossible. Thus, $\{x_a\}$ is bounded uniformly for $a\nearrow a^*$.
Similarly, one can check that $\lim_{a\nearrow a^*}V_{\Omega}(x_a)=0$. This completes the proof of Theorem \ref{thm:1.2}.
\qed

\begin{appendix}

\section{Appendix}

In this appendix, we
give the detailed proof of some results used in the proof of Lemma \ref{pro:3.2} for the reader's convenience.

\begin{lem}\label{lem:A.3}
Let $\varepsilon_a>0$ and $w_{a}(x)$ be defined by \eqref{3.1} and \eqref{3.2}, respectively. Then $\varepsilon_a^2\int_{\R^2}\ln|x-y||w_a(y)|^2dy\in C^{\alpha}_{loc}(\R^2)$ and
$\big\{\varepsilon_a^2\int_{\R^2}\ln|x-y||w_a(y)|^2dy\big\}$ is bounded uniformly in $C^{ \alpha}_{loc}(\R^2)$ as $a\nearrow a^*$.
\end{lem}
\vskip 0.1truein
\noindent\textbf{Proof.} 
Firstly, for any $R>0$, note from \cite[Theorem 4.2]{Gilbarg} that
\begin{equation}\label{a7}
\begin{split}
g_{a}(x):=\varepsilon_a^2\int_{B_{R+1}(0)}\ln|x-y||w_a(y)|^2dy\in C^{2}\big(B_{R+1}(0)\big),
\end{split}
\end{equation}
and $g_{a}(x)$ satisfies
\begin{equation}\label{a8}
\begin{split}
\Delta g_{a}(x)=2\pi\varepsilon_a^2|w_a(x)|^2\ \ \text{in}\ \ B_{R+1}(0).
\end{split}
\end{equation}
Moreover, we derive from (\ref{3.451}) that
\begin{equation}\label{a9}
\begin{split}
|g_{a}(x)|=&\Big|\varepsilon_a^2\int_{B_{R+1}(0)}\ln|x-y||w_a(y)|^2dy\Big|\\
\leq&\int_{B_{R+1}(0)}\ln\Big(1+|x-y|\Big)|w_a(y)|^2dy\\
&+\int_{B_{R+1}(0)}\ln\Big(1+\frac{1}{|x-y|}\Big)|w_a(y)|^2dy\\
\leq&\int_{\R^2}\ln\Big(1+|x-y|\Big)|w_a(y)|^2dy
+\int_{\R^2}\ln\Big(1+\frac{1}{|x-y|}\Big)|w_a(y)|^2dy\\
\leq& a^*\ln\big(2+R\big)+C\ \ \text{in}\ \ B_{R+1}(0)\ \  \text{as}\ \ a\nearrow a^*.
\end{split}
\end{equation}
Note that $\big\{\varepsilon_a^2|w_a(x)|^2\big\}$ is bounded uniformly in $C^{\alpha}(B_{R+1}(0))$ as $a\nearrow a^*$, it then follows from (\ref{a9}) and the Schauder estimate \cite[Theorem 6.2]{Gilbarg} that $\big\{g_{a}(x)\big\}$ is bounded uniformly in $C^{2,\alpha}(B_{R}(0))$ as $a\nearrow a^*$.

On the other hand, for any $x_1, x_2\in B_{R}(0)$, direct calculations yield that
\begin{equation}\label{a10}
\begin{split}
&\frac{\Big|\varepsilon_a^2\int_{\mathbb{R}^{2}\setminus B_{R+1}(0)}\ln|x_{2}-y||w_a(y)|^2dy
-\varepsilon_a^2\int_{\mathbb{R}^{2}\setminus B_{R+1}(0)}\ln|x_{1}-y||w_a(y)|^2dy\Big|}{\big|x_2-x_1\big|^{\alpha}}\\
\leq&|x_2-x_1|^{-\alpha}\int_{\mathbb{R}^{2}\setminus B_{R+1}(0)}\big|\ln|x_{2}-y|-\ln|x_{1}-y|\big||w_a(y)|^2dy\\
\leq&|x_2-x_1|^{1-\alpha}\int_{\mathbb{R}^{2}\setminus B_{R+1}(0)}\frac{|w_a(y)|^2}
{\big|x_1+\theta(x_2-x_1)-y\big|}dy\\
\leq&|x_2-x_1|^{1-\alpha}\int_{\mathbb{R}^{2}\setminus B_{R+1}(0)}|w_a(y)|^2dy\\
\leq&C|x_2-x_1|^{1-\alpha}\leq C(R)\ \ \text{in}\ \ B_{R}(0)\ \  \text{as}\ \ a\nearrow a^*,
\end{split}
\end{equation}
where $\theta\in(0,1)$. Moreover, similar to (\ref{a9}), we can obtain that
\begin{equation}\label{a11}
\begin{split}
\Big|\varepsilon_a^2\int_{\mathbb{R}^{2}\setminus B_{R+1}(0)}\ln|x-y||w_a(y)|^2dy\Big|
\leq a^*\ln\big(1+R\big)+C\ \ \text{in}\ \ B_{R}(0)\ \  \text{as}\ \ a\nearrow a^*.
\end{split}
\end{equation}
Thus, we deduce from (\ref{a10}) and (\ref{a11}) that $$\varepsilon_a^2\int_{\mathbb{R}^{2}\setminus B_{R+1}(0)}\ln|x-y||w_a(y)|^2dy\in C^{\alpha}(B_{R}(0))$$
 and $\{\varepsilon_a^2\int_{\mathbb{R}^{2}\setminus B_{R+1}(0)}\ln|x-y||w_a(y)|^2dy\}$ is bounded uniformly in $C^{\alpha}(B_{R}(0))$ as $a\nearrow a^*$.

Since $R>0$ is arbitrary, we get from the above facts that $$\varepsilon_a^2\int_{\R^2}\ln|x-y||w_a(y)|^2dy\in C^{\alpha}_{loc}(\R^2)$$
and
$\big\{\varepsilon_a^2\int_{\R^2}\ln|x-y||w_a(y)|^2dy\big\}$ is bounded uniformly in $C^{ \alpha}_{loc}(\R^2)$ as $a\nearrow a^*$. This completes the proof of Lemma \ref{lem:A.3}.
\qed

\end{appendix}

\vspace {.5cm}
\noindent {\bf Acknowledgements:}  The authors thank Professor Yujin Guo very much for his fruitful discussions on the present paper.

\end{document}